\def\ppp{}
\def\pp{}
\def\section{}

\def\circv{\breve v}

\def\siml#1{{\overset{L^{\#}}\to\sim}}

\def\re{\operatorname{Re}}

\def\ad{\operatorname{ad}}

\def\tr{\operatorname{tr}}

\def\normi#1{\|#1\|_{\infty}}

\def\pmtwo#1#2#3#4{\pmatrix#1&#2\\#3&#4\endpmatrix}

\def\cb{{\Bbb C}}
\def\rb{{\Bbb R}}





\def\by{{\bold y}}

\let\tildesymbol=\~

\magnification=1000
\input amstex
\input pictex

\documentstyle{amsppt}
\loadbold
\baselineskip = 18pt
\overfullrule = 0pt

\vsize=9.0truein
\hsize=6.5truein

\define\n{\noindent}
\define\bb#1{\Bbb{#1}}
\define\cl#1{\Cal{#1}}

\define\ovl{\overline}

\define\vp{\varepsilon}
\define\sst{\scriptstyle}

\def\mnorm#1{\langle#1\rangle_\infty}

\rightheadtext{A priori $L^p$ estimates for solutions of
Riemann--Hilbert problems} \topmatter
\title
A priori ${\bold L}^{\bold p}$ estimates for solutions of
Riemann--Hilbert Problems
\endtitle
\author
Percy Deift and Xin Zhou
\endauthor
\affil
New York University \quad Duke University
\endaffil
\thanks
The work of the first author was supported in part by NSF Grants DMS~0003268 and DMS~0296084 and the work of the second author was supported in part by NSF Grant DMS~0071398.
\endthanks
\endtopmatter

\document
\def\pp{}
\def\ppp{}

\head 1. Introduction\endhead In this paper we prove a general
result establishing {\it a priori\/} $L^p$ estimates for solutions
of Riemann-Hilbert Problems (RHP's) in terms of auxiliary
information involving an associated ``conjugate" problem (see
Conjugation Lemma \pp{1.39} below). We then use the result to
obtain uniform estimates for a RHP (see Theorem \pp{1.48}) that
plays a crucial role in analyzing the long-time behavior of
solutions of the perturbed nonlinear Schr\"odinger equation on the
line. Theorem \pp{1.48} is proved by combining Conjugation Lemma
\pp{1.39} with the steepest-descent method for RHP's introduced by
the authors in [DZ1]. We do not apply the steepest-descent method
directly to Theorem \pp{1.48}. Rather, as explained in the text,
we proceed by rephrasing Theorem \pp{1.48} as an equivalent
inhomogeneous RHP in which the underlying objects $M_\pm$ (see
Theorem \pp{1.52}) have appropriate analyticity properties and can
be deformed around the stationary phase point much as in the
manner of the classical method of
stationary-phase/steepest-descent.

We begin by introducing a variety of definitions and results that
arise in the theory of Riemann-Hilbert Problems (RHP's).

 Let $\Sigma$ be an oriented contour in $\bb C$ and consider the associated Cauchy operator
$$Ch(z) = C_\Sigma h(z) \equiv \int_\Sigma \frac{h(s)}{s-z} \frac{ds}{2\pi i},\qquad z\in \bb C\backslash\Sigma.$$
If we move along the contour in the direction of the orientation,
we say, by convention,  that the $(+)$-side (resp. $(-)$-side)
lies to the left (resp.\ right). The following properties and
estimates which will be used without further comment throughout
the text are true for a very general class of contours (see, for
example,[DS],[Dur]). This class certainly includes contours that
are finite unions of smooth curves in $\ovl{\bb C}$ such that
$\ovl\cb\setminus\Sigma$ has a finite number of components as in
Figures \pp{1.35}, \pp{1.38} and \pp{3.35} below.
\medskip

\n (\ppp{1.1})\ Let $h\in L^p(\Sigma, |dz|), 1\le p<\infty$. Then
$$C^\pm h(z) \equiv \lim_{\sst z'\to z\atop \sst z'\in (\pm)-\text{side of } \Sigma} (Ch)(z')$$
exists as a non-tangential limit for a.e.\ $z\in\Sigma$.

\n{(\ppp{1.2})} Let $h\in L^p(\Sigma, |dz|), 1<p<\infty$. Then

$$\|C^\pm h\|_{L^p(\Sigma)} \le c_p\|h\|_{L^p(\Sigma)}$$
for some constant $c_p$.

\n{(\ppp{1.3})} $C^\pm=\pm1/2-H/2$, where $H$ is the Hilbert
transform
$$Hf(z)=P.V.\int_\Sigma\frac{f(s)}{z-s}\frac{ds}{i\pi}.$$

By (\pp{1.2}), for $h\in L^p(\Sigma,|dz|)$, $1<p<\infty$,
$$\|Hh\|_{L^p(\Sigma)}\le c_p\|h\|_{L^p(\Sigma)}$$ for some
constant $c_p$.

\n{(\ppp{1.4})} $C^+-C^- = 1,\ 1\le p < \infty$.\medskip

\n Let $v$ be a $k\times k$ {\it jump matrix\/} on $\Sigma$ i.e.\ $v$ is a measurable map from $\Sigma \to GL(k,\bb C)$ with $v,v^{-1}\in L^\infty(\Sigma\to GL(k,\bb C))$. Define the associated singular integral operator

$$C_vh = C^-(h(v-I))\tag{\ppp{1.5}}$$
acting on $L^p(\Sigma)$-matrix-valued functions. Clearly $C_v$ is
bounded from $L^p\to L^p$ for all $1<p<\infty$.

Given $\Sigma$ and $v$, the operator $C_v$ is intimately connected
with the solution of associated RHP's on $\Sigma$. For
$1<p<\infty$, we say that a pair of $L^p(\Sigma)$-functions $f_\pm
\in \partial C(L^p)$ if there exists a (unique) function $h\in
L^p(\Sigma)$ such that

$$f_\pm(z) = (C^\pm h)(z),\qquad z\in \Sigma.\tag\ppp{1.6}$$
In turn we call $f(z) = Ch(z), z\in \bb C\backslash\Sigma$, the
{\it extension\/} of $f_\pm = C^\pm h\in\partial C(L^p)$ off
$\Sigma$. Observe that if $F(z)=(Cf)(z)$ for $f\in L^p(\Sigma)$
and $G(z)=(Cg)(z)$ for $g\in L^q(\Sigma)$, $1<p,q<\infty$, then a
simple computation shows that if $1/r=1/p+1/q\le1$,
$$FG(z)=(Ch)(z),$$where
$$h=-\frac12((Hf)g+fHg)\in L^r(\Sigma).\tag\ppp{1.7}$$
If $F=Cf(z)$ for $f\in L^p(\Sigma)$, $1<p<\infty$, and $G$ is
analytic and bounded in $\cb\setminus\Sigma$, then
$$FG(z)=Ch(z),\tag\ppp{1.8}$$
where
$$ h=(C^+f)G_+-(C^-f)G_-\in L^p(\Sigma)\tag\ppp{1.9}$$
and $G_\pm$ denote the non-tangential limits (see [Dur]) of $G$ on
$\Sigma$.

\definition{Definition \ppp{1.10}} {\bf Inhomogeneous Riemann--Hilbert Problem IRHP$_{\bold{L^p}}$}
(see [CG], [DZ5], [DZ7]). Fix $1<p<\infty$. Given $\Sigma,v$ and a
function $F\in L^p(\Sigma)$, we say that $M_\pm \in \partial
C(L^p)$ solves an IRHP$_{L^p}$ if

$$M_+(z) = M_-(z) v(z) + F(z),\qquad z\in \Sigma.\tag\ppp{1.11}$$

If $M(z)$ is the extension of $M_\pm$ off $\Sigma$, we see that

\item{$\bullet$} $M(z)$ is analytic in $\bb C\backslash \Sigma$,

\item{$\bullet$} $M_+(z) = M_-(z) v(z) + F(z), z\in\Sigma$, where
$M_\pm(z) = \lim\limits_{\sst z'\to z\atop z'\in (\pm)\text{-side
of } \Sigma} M(z'),$

\item{$\bullet$} $M(z)\to 0$ as $z\to\infty$ in any non-tangential
direction.

\enddefinition

\remark{Remark}

In [DZ4] and [DZ5], IRHP$_{L^p}$ is referred to as an
inhomogeneous RHP of type 2, IRHP2$_{L^p}$. There is also the
notion of an inhomogeneous RHP of type 1, IRHP1$_{L^p}$. As we
only need type 2 in this paper, we simply drop the ``2''.

\endremark

The following relations are basic.

\proclaim{Proposition \ppp{1.12}}

Let $1<p<\infty$. Then $1-C_v$ is a bijection in $L^p(\Sigma)$ if
and only if IRHP$_{L^p}$ has a unique solution for all $F\in
L^p(\Sigma)$. Moreover if $(1-C_v)^{-1}$ exists, then for $h\in
L^p(\Sigma)$

$$(1-C_v)^{-1} h = (M_++h)v^{-1} = (M_-+h)\tag\ppp{1.13}$$
where $M_\pm$ solves IRHP$_{L^p}$ with $F = h(v-I)$. Conversely,
if $M_\pm$ solves IRHP$_{L^p}$ with $F\in L^p(\Sigma)$, then

$$M_+ = ((1-C_v)^{-1} (C^-F)) v + F \quad \text{and}\quad M_- = (1-C_v)^{-1} C^-F.\tag\ppp{1.14}$$

\endproclaim

Suppose $v = (v^-)^{-1}v^+$ is a factorization of $v$ with
$v^\pm$, $(v^\pm)^{-1}\in L^\infty(\Sigma, GL(k,\bb C))$. Set
$w^+\equiv v^+-I$, $w^-\equiv I-v^-$ and let $w = (w^-, w^+)$.
Define the associated singular integral operator (cf.\ (\pp{1.5}))

$$C_wh = C^+(hw^-) + C^-(hw^+)\tag\ppp{1.15}$$
acting on $L^p(\Sigma)$-matrix-valued function $h$. Note that for
the trivial factorization $v = (I)^{-1}v$, $w = (0,v-I)$, we have
$C_w=C_v$.

The following result describes the relation between $1-C_w$ and
$1-C_{w'}$ for two different factorizations $v = (v^-)^{-1} v^+ =
(v^{\prime-})^{-1} v^{\prime +}$ of $v$.

\proclaim{Proposition \ppp{1.16} (see [Z1],[DZ4], [DZ5])}

Suppose $1<p<\infty$. The operator $1-C_w$ is bijective in
$L^p(\Sigma)$ for all factorization $v = (v^-)^{-1}v^+ =
(I-w^-)^{-1}(I+w^+)$, if and only if $(1-C_{w'})$ is bijective for
at least one factorization $v  =(v^{\prime -})^{-1} v^{\prime +} =
(I-w^{\prime -})^{-1}(I+w^{\prime +})$. Moreover, for $f\in
L^p(\Sigma)$

$$(1-C_w)^{-1} f = ((1-C_{w'})^{-1}f)b$$
where $b = v^{\prime +}(v^+)^{-1} = v^{\prime-}(v^-)^{-1}$.

\endproclaim

\subheading{Duality}A non--degenerate bilinear pairing for vector
functions in $L^p(\Sigma)$ and $L^q(\Sigma)$, $1<p,q<\infty$,
$1/p+1/q=1$, is given by
$$\langle f,g\rangle=\int_\Sigma\tr(f(z)g(z)^T)dz.\tag\ppp{1.17}$$

Using (\pp{1.3}), we see that with respect to this pairing, the
dual operators $(C^\pm)'$ are given by
$$(C^\pm)'=-C^\mp\tag\ppp{1.18}$$
and
$$C_v'h=(C^+h)(I-v^T)=R_{I-v^T}C^+h,\tag\ppp{1.19}$$
where $R_A$ denotes multiplication on the right by a matrix $A$.
Now $1-C_v$ is bijective in $L^p(\Sigma)$ if and only if $1-C_v'$
is bijective in $L^q(\Sigma)$, $1/p+1/q=1$, $1<p,q<\infty$, and
$$\|(1-C_v')^{-1}\|_{L^q\to L^q}=\|(1-C_v)^{-1}\|_{L^p\to L^p}.\tag\ppp{1.20}$$
On the other hand, it is well-known that if $K$ and $L$ are
bounded linear operators in a Banach space, then (see e.g. [D] and
the references therein) $1-KL$ is invertible if and only if $1-LK$
is invertible and $(1-KL)^{-1}=1+K (1-LK)^{-1}L$. Setting $K=C^+$
and $L=R_{I-v^T}$, we find that $1-C^+R_{I-v^T}$ is a bijection in
$L^q(\Sigma)$ and
$$(1-C^+R_{I-v^T})^{-1}=1+C^+(1-C_v')^{-1}R_{I-v^T}\tag\ppp{1.21}$$
Observe that $(C^+R_{I-v^T})h=C^+h(I-v^T)$, which is precisely an
operator of the form $C_w$ in (\pp{1.15}) above corresponding to
the jump matrix $(v^T)^{-1}=(v^-)^{-1}v^+$ with
$w=(w^-,w^+)=(I-v^T,0)$.

We need the following simple facts. If $|\cdot|$ denotes any
sub-multiplicative norm on an algebra, $|AB|\le|A||B|$, then
$|Id|\ge1$ and if A is invertible,
$$|A||A^{-1}|\ge1\tag\ppp{1.22}$$
and hence
$$\mnorm v\equiv \max(\normi v,\normi{v^{-1}})\ge1.\tag\ppp{1.23}$$
If $1-C_v$ is invertible in $L^p(\Sigma)$, then a simple
calculation using the identity $(1-C_v)^{-1}(1-C_v)=1$, together
with (\pp{1.23}), shows that
$$0<c\le\|(1-C_v)^{-1}\|_{L^p(\Sigma)\to L^p(\Sigma)}\mnorm v\tag\ppp{1.24}$$
for some constant $c=c_p$ independent of $v$.

Throughout the text constants $c$ will be used generically.
Statements such as $\|f\| \le 2c (1+e^c) \le c$, for example,
should not cause any confusion.

 \subheading{Normalized RHP} Let $\Sigma$ be an oriented contour
with associated jump matrix $v$. Suppose $v-I\in L^p(\Sigma)$ for
some $1<p<\infty$. We say that $\phi_\pm$ solves the {\it
normalized RHP\/} $(\Sigma,v)_{L^p}$ if $\phi_\pm -I\in \partial
C(L^p)$ solves the IRHP$_{L^p}$

$$\phi_+-I = (\phi_- - I)v + (v-I) \quad \text{in}\quad \Sigma\tag\ppp{1.25}$$
or equivalently

$$\phi_+ = \phi_-v \quad \text{on}\quad \Sigma.\tag{\ppp{1.26}}$$

If $\phi(z)-I$ is the extension of $\phi_\pm-I$ off $\Sigma$, we
see that

\item{$\bullet$} $\phi(z)$ is analytic in $\bb C\backslash
\Sigma$,

\item{$\bullet$} $\phi_+(z) = \phi_-(z) v(z), z\in \Sigma$,
$\text{where} \quad \phi_\pm(z)  = \lim_{\sst z'\to z \atop \sst
z'\in (\pm)\text{-side of } \Sigma} \phi(z').$

\item{$\bullet$} $\phi(z)\to I$ as $z\to \infty$ in any
non-tangential direction.

If $1-C_v$ is invertible in $L^p$ then, by Proposition \pp{1.12},
the normalized solution exists, is unique, and is given by
$$\phi_-=I+(1-C_v)^{-1}C^-(v-I)\tag\ppp{1.27}$$
and hence $\|\phi_--I\|_{L^p}\le c\|(1-C_v)^{-1}\|_{L^p\to
L^p}\|v-I\|_{L^p}$. Using the relation (\pp{1.26}) together with
(\pp{1.23},\pp{1.24}), we then obtain the inequalities
$$\|\phi_\pm-I\|_{L^p}\le c\|(1-C_v)^{-1}\|_{L^p\to L^p}\mnorm v
\|v-I\|_{L^p}.\tag\ppp{1.28}$$ As noted above, the operator
$1-C^+R_{I-v^T}$ associated with the jump matrix $(v^T)^{-1}$ is
invertible in $L^q(\Sigma)$. Let $\psi_\pm$ solve the normalized
RHP $(\Sigma,v)_{L^q}$, $\psi_+=\psi_-(v^T)^{-1}$, $\psi_\pm-I\in
\partial C(L^q)$. As in (\pp{1.27}), we have
$\psi_\pm=I+(1-C^+R_{I-v^T})^{-1}C^+(I-v^T)$, which leads as above
to the bound
$$\align\|\psi_+-I\|_{L^q}&\le c\|(1-C^+R_{I-v^T})^{-1}\|_{L^q\to
L^q}\|v-I\|_{L^q}\tag\ppp{1.29}\\
&\le c\|(1-C'_v)^{-1}\|_{L^q\to L^q} \mnorm v\|v-I\|_{L^q},\ \
\text{ by (\pp{1.21})},\\
&=c\|(1-C_v)^{-1}\|_{L^p\to L^p} \mnorm v\|v-I\|_{L^q},\ \ \text{
by (\pp{1.20})},
\endalign$$
which leads in turn to the estimate $$ \|\psi_--I\|_{L^q}\le
c\|(1-C_v)^{-1}\|_{L^p\to L^p} \mnorm
v^2\|v-I\|_{L^q}\tag\ppp{1.30}$$ as before. Using the jump
relations for $\phi_\pm$ and $\psi_\pm$, we see that
$\psi_+\phi_+^T=\phi_-vv^{-1}\psi_-^T=\phi_-\psi_-^T$. But
$\phi_\pm\psi_\pm^T-I=C^\pm h$, where $h\in L^1+L^p+L^q$, by
(\pp{1.7). Hence $h=C^+h-C^-h=0$, and so
$\psi_\pm^T=\phi_\pm^{-1}$. We conclude that
$$\|\phi_+^{-1}-I\|_{L^q}\le c \|(1-C_v)^{-1}\|_{L^p\to L^p}\mnorm
v\|v-I\|_{L^q}\tag\ppp{1.31}$$ and
$$\|\phi_-^{-1}-I\|_{L^q}\le
c \|(1-C_v)^{-1}\|_{L^p\to L^p}\mnorm
v^2\|v-I\|_{L^q}.\tag\ppp{1.32}$$

Recall that a contour $\Gamma\subset \bb C$ is {\it complete\/}
(see e.g.\ [Z2]) if

$$\align
&\bb C\backslash \Gamma \text{ is a disjoint union of two,
possibly disconnected,
open regions $\Omega_+$ and $\Omega_-$, and}\tag\ppp{1.33}\\
&\Gamma \text{ may be viewed as the positively oriented boundary for $\Omega_+$ and also as the}\tag\ppp{1.34}\\
&\text{negatively oriented boundary of $\Omega_-$.}
\endalign$$

\n Two examples of such contours are $\Gamma=\bb R$ and $\Gamma =
\bb R\cup i\bb R$.

\bigskip
\font\thinlinefont=cmr5

\centerline{\beginpicture \setcoordinatesystem units <.7cm,.7cm>
\linethickness=1pt
\setshadesymbol ({\thinlinefont .}) \setlinear
%
%
\linethickness= 0.500pt \setplotsymbol ({\thinlinefont .})
\putrule from  2.540 21.590 to  7.620 21.590
%
%
\linethickness= 0.500pt \setplotsymbol ({\thinlinefont .})
\putrule from 12.700 21.590 to 17.780 21.590
%
%
\linethickness= 0.500pt \setplotsymbol ({\thinlinefont .})
\putrule from 15.240 24.130 to 15.240 19.050
%
%
\linethickness= 0.500pt \setplotsymbol ({\thinlinefont .})
\putrule from  4.286 21.590 to  4.604 21.590
%
%
\plot  4.350 21.526  4.604 21.590  4.350 21.654 /
%
%
%
\linethickness= 0.500pt \setplotsymbol ({\thinlinefont .})
\putrule from 16.034 21.590 to 16.192 21.590
%
%
\plot 15.939 21.526 16.192 21.590 15.939 21.654 /
%
%
%
\linethickness= 0.500pt \setplotsymbol ({\thinlinefont .})
\putrule from 14.129 21.590 to 13.970 21.590
%
%
\plot 14.224 21.654 13.970 21.590 14.224 21.526 /
%
%
%
\linethickness= 0.500pt \setplotsymbol ({\thinlinefont .})
\putrule from 15.240 23.178 to 15.240 22.860
%
%
\plot 15.176 23.114 15.240 22.860 15.304 23.114 /
%
%
%
\linethickness= 0.500pt \setplotsymbol ({\thinlinefont .})
\putrule from 15.240 20.003 to 15.240 20.161
%
%
\plot 15.304 19.907 15.240 20.161 15.176 19.907 /
%
%
%
\put{$\Omega_+$} [lB] at  4.128 22.543
%
%
\put{$\Omega_-$} [lB] at  4.128 20.320
%
%
\put{$\Omega_-$} [lB] at 13.811 22.543
%
%
\put{$\Omega_+$} [lB] at 16.034 22.543
%
%
\put{$\Omega_+$} [lB] at 13.811 20.320
%
%
\put{$\Omega_-$} [lB] at 16.034 20.320 \linethickness=0pt
\putrectangle corners at  2.515 24.155 and 17.805 19.025
\endpicture}

\bigskip

\centerline{Figure \ppp{1.35} $\Gamma = \bb R,\ \Gamma = \bb R \cup
i\bb R$}\bigskip

Note that a complete contour $\Gamma$ comes equipped with two
natural orientations (we may always relabel $\Omega_\pm$ as
$\Omega_\mp$). Unless stated otherwise, we will always choose one
of these orientations for the specification of RHP's on $\Gamma$.
Note also that if $\Gamma$ is complete, then by Cauchy's theorem
$C^+C^- = C^-C^+  = 0$, and hence by the basic relation
$C^+-C^-=1$, we see that $C^+$ and $-C^-$ are complementary
projections in $L^p(\Gamma)$, $1<p<\infty$.

In this section, we will also consider extended contours
$\Gamma_{\text{ext}} = \Gamma \cup \Gamma'$ where

$$\align
&\text{dist}(\Gamma,\Gamma')>0\tag\ppp{1.36}\\
&\Gamma_{\text{ext}} \text{ is complete.}\tag\ppp{1.37}
\endalign$$

Two examples of such extended contours are extensions of $\bb R$
and $\bb R\cup i\bb R$ as in Figure \pp{1.38}.

\vskip2in

\bigskip\bigskip
\font\thinlinefont=cmr5 \centerline{\beginpicture
\setcoordinatesystem units <.7cm,.7cm>
\linethickness=1pt
\setshadesymbol ({\thinlinefont .}) \setlinear
%
%
\linethickness= 0.500pt \setplotsymbol ({\thinlinefont .})
\putrule from  2.540 21.590 to  7.620 21.590
%
%
\linethickness= 0.500pt \setplotsymbol ({\thinlinefont .})
\putrule from 12.700 21.590 to 17.780 21.590
%
%
\linethickness= 0.500pt \setplotsymbol ({\thinlinefont .})
\putrule from 15.240 24.130 to 15.240 19.050
%
%
\linethickness= 0.500pt \setplotsymbol ({\thinlinefont .})
\putrule from  4.286 21.590 to  4.604 21.590
%
%
\plot  4.350 21.526  4.604 21.590  4.350 21.654 /
%
%
%
\linethickness= 0.500pt \setplotsymbol ({\thinlinefont .})
\putrule from 16.034 21.590 to 16.192 21.590
%
%
\plot 15.939 21.526 16.192 21.590 15.939 21.654 /
%
%
%
\linethickness= 0.500pt \setplotsymbol ({\thinlinefont .})
\putrule from 14.129 21.590 to 13.970 21.590
%
%
\plot 14.224 21.654 13.970 21.590 14.224 21.526 /
%
%
%
\linethickness= 0.500pt \setplotsymbol ({\thinlinefont .})
\putrule from 15.240 23.178 to 15.240 22.860
%
%
\plot 15.176 23.114 15.240 22.860 15.304 23.114 /
%
%
%
\linethickness= 0.500pt \setplotsymbol ({\thinlinefont .})
\putrule from 15.240 20.003 to 15.240 20.161
%
%
\plot 15.304 19.907 15.240 20.161 15.176 19.907 /
%
%
%
\linethickness= 0.500pt \setplotsymbol ({\thinlinefont .})
\putrule from  2.540 24.130 to  7.620 24.130
%
%
\linethickness= 0.500pt \setplotsymbol ({\thinlinefont .})
\putrule from  2.540 19.050 to  7.620 19.050
%
%
\linethickness= 0.500pt \setplotsymbol ({\thinlinefont .})
\putrule from  4.604 24.130 to  4.445 24.130
%
%
\plot  4.699 24.194  4.445 24.130  4.699 24.066 /
%
%
%
\linethickness= 0.500pt \setplotsymbol ({\thinlinefont .}) \plot
4.921 19.050  4.921 19.050 /
%
%
\linethickness= 0.500pt \setplotsymbol ({\thinlinefont .})
\putrule from  5.080 19.050 to  4.763 19.050
%
%
\plot  5.017 19.114  4.763 19.050  5.017 18.986 /
%
%
%
\linethickness= 0.500pt \setplotsymbol ({\thinlinefont .})
\putrule from 12.700 22.860 to 13.970 22.860 \putrule from 13.970
22.860 to 13.970 24.130
%
%
\linethickness= 0.500pt \setplotsymbol ({\thinlinefont .})
\putrule from 16.510 24.130 to 16.510 22.860 \putrule from 16.510
22.860 to 17.780 22.860
%
%
\linethickness= 0.500pt \setplotsymbol ({\thinlinefont .})
\putrule from 12.700 20.320 to 13.970 20.320 \putrule from 13.970
20.320 to 13.970 19.050
%
%
\linethickness= 0.500pt \setplotsymbol ({\thinlinefont .})
\putrule from 16.510 19.050 to 16.510 20.320 \putrule from 16.510
20.320 to 17.780 20.320
%
%
\linethickness= 0.500pt \setplotsymbol ({\thinlinefont .})
\putrule from 13.018 22.860 to 13.176 22.860
%
%
\plot 12.922 22.796 13.176 22.860 12.922 22.924 /
%
%
%
\linethickness= 0.500pt \setplotsymbol ({\thinlinefont .})
\putrule from 17.462 22.860 to 17.145 22.860
%
%
\plot 17.399 22.924 17.145 22.860 17.399 22.796 /
%
%
%
\linethickness= 0.500pt \setplotsymbol ({\thinlinefont .})
\putrule from 12.859 20.320 to 13.176 20.320
%
%
\plot 12.922 20.256 13.176 20.320 12.922 20.384 /
%
%
%
\linethickness= 0.500pt \setplotsymbol ({\thinlinefont .})
\putrule from 17.462 20.320 to 17.304 20.320
%
%
\plot 17.558 20.384 17.304 20.320 17.558 20.256 /
%
%
%
\put{$\Omega_+$} [lB] at  4.128 22.543
%
%
\put{$\Omega_-$} [lB] at  4.128 20.320
%
%
\put{$\Omega_-$} [lB] at  4.128 24.448
%
%
\put{$\Omega_+$} [lB] at  3.969 18.098
%
%
\put{$\Gamma'$} [lB] at  1.587 23.971
%
%
\put{$\Gamma$} [lB] at  1.587 21.431
%
%
\put{$\Gamma'$} [lB] at  1.587 18.891
%
%
\put{$\Gamma'$} [lB] at 11.748 22.701
%
%
\put{$\Gamma'$} [lB] at 11.748 20.161
%
%
\put{$\Gamma'$} [lB] at 18.415 22.701
%
%
\put{$\Gamma$} [lB] at 18.415 21.431
%
%
\put{$\Gamma'$} [lB] at 18.574 20.161
%
%
\put{$\Omega_+$} [lB] at 13.018 23.336
%
%
\put{$\Omega_-$} [lB] at 16.828 23.336
%
%
\put{$\Omega_-$} [lB] at 12.859 19.526
%
%
\put{$\Omega_+$} [lB] at 16.828 19.526
%
%
\put{$\Omega_-$} [lB] at 14.287 22.225
%
%
\put{$\Omega_+$} [lB] at 15.716 22.225
%
%
\put{$\Omega_+$} [lB] at 14.287 20.637
%
%
\put{$\Omega_-$} [lB] at 15.716 20.637 \linethickness=0pt
\putrectangle corners at  1.587 24.733 and 18.574 18.098
\endpicture}

\centerline{Figure \ppp{1.38} $\Gamma_{\text{ext}} = \Gamma\cup
\Gamma'$}\bigskip

Given $\Gamma_{\text{ext}} = \Gamma\cup\Gamma'$, let $R$ and
$R^{-1}$ be analytic and bounded matrix functions in $\bb
C\backslash \Gamma_{\text{ext}}$. As noted above following
(\pp{1.9}), $R$ and $R^{-1}$ have non-tangential boundary values
almost everywhere on $\Gamma_{\text{ext}}$. We denote the boundary
values from $\Omega_+$ by $R_+,R^{-1}_+$ and from $\Omega_-$ by
$R_-$, $R^{-1}_-$. The following Lemma provides bounds on
$\|(1-C_v)^{-1}\|_{L^p(\Gamma)\to L^p(\Gamma)}$ in terms of
auxiliary information involving an associated ``conjugate"
problem, and is the main structural result in this paper.

\proclaim{Conjugation Lemma \ppp{1.39}}

Let $2<p<\infty$ and let $\Gamma_{\text{ext}} = \Gamma \cup
\Gamma'$, $R,R^{-1}$ be given as above. Let $v,\breve v\colon \
\Gamma\to GL(k,\bb C)$ be jump matrices that are related via

$$v  = R^{-1}_-\breve vR_+.\tag\ppp{1.40}$$
Assume in addition that

$$v-I, \ \breve v-I \in L^2(\Gamma).\tag\ppp{1.41}$$
Suppose that

$$(1-C_{\breve v})^{-1} \text{ exists in $L^2(\Gamma)$ and also in $L^p(\Gamma)$}\tag\ppp{1.42}$$
and that

$$(1-C_v)^{-1} \text{ exists in } L^2(\Gamma).\tag\ppp{1.43}$$
Then $(1-C_v)^{-1}$ exists in $L^p(\Gamma)$ and

$$\|(1-C_v)^{-1}\|_{L^p(\Gamma)\to L^p(\Gamma)} \le c^{\#}_v,\tag\ppp{1.44}$$
where

$$\align
c^\#_v &= c\|R\|_{L^\infty(\bb C\backslash \Gamma_{\text{ext}})}
\|R^{-1}\|_{L^\infty(\bb C\backslash \Gamma_{\text{ext}})}
\|(1-C_{\breve v})^{-1}\|_{L^p (\Gamma)\to L^p(\Gamma)}
\|(1-C_{\breve v})^{-1}\|_{L^2 (\Gamma)\to L^2(\Gamma)}\tag\ppp{1.45}\\
&\quad \times \|(1-C_v)^{-1}\|^2_{L^2(\Gamma)\to L^2(\Gamma)}
\mnorm v^3 \mnorm {\breve v}^2  (1+ \|\breve
v-I\|_{L^2(\Gamma)})^2  (1+\|v-I\|_{L^2(\Gamma)})^2.
\endalign$$

\endproclaim

\remark{Remark \ppp{1.46}}

\noindent(a) \ The assumption that $\Gamma$ and
$\Gamma_{\text{ext}}$ are complete is not necessary. However, if
$\Gamma_{\text{ext}}$ is not complete then the meaning of $R_\pm$
is not clear and the precise statement of the theorem becomes more
complex. Of course, any contour can be extended trivially to a
complete contour.

\noindent(b) \ If $v,\breve v$ satisfy all the assumptions in the
Conjugation Lemma \pp{1.39}, except $\breve v$ is replaced by
$(\breve v^T)^{-1}$ in condition (\pp{1.42}), then the Conjugation
Lemma together with the  duality relations given above implies
that $(1-C_v)^{-1}$ exists in $L^q(\Gamma)$, $1/p+1/q=1$, and
$$\|(1-C_v)^{-1}\|_{L^q(\Gamma)\to L^q(\Gamma)}\le c^\#_v,$$ where $c^\#_v$ has
similar structure to the constant in (\pp{1.45}). We leave the
details to the reader.

\endremark

An application of the Conjugation Lemma to PDE's is given in
Theorem
 \pp{1.48} below, as noted at the beginning of this section.

 Let $H^1$ denote the first Sobolev space $\{f\colon\ f,f'\in
L^2(\bb R)\}$ and define $H^1_1 \equiv H^1\cap \{f\colon \
\|f\|_{L^\infty(\bb R)} <1\}$. For $\Sigma=\bb R$, oriented from
$-\infty\to\infty$, consider the $2\times 2$ jump matrix

$$v_\theta = \left(\matrix
1-|r(z)|^2&r(z)e^{i\theta}\\
-\ovl{r(z)} e^{-i\theta}&1\endmatrix\right),\qquad z\in \bb
R,\tag{\ppp{1.47}}$$ where $r\in H^1_1$ and $\theta = xz-tz^ 2$,
$x,t\in \bb R$.

\proclaim{Theorem \ppp{1.48}}

Suppose $r\in H^1_1, \|r\|_{H^1}\le \lambda, \|r\|_{L^\infty} \le
\rho <1$. Then for $x,t\in \bb R$, and for any $2\le p < \infty$,
$(1-C_{v_\theta})^{-1}$ exists as a bounded operator in $L^p(\bb
R)$, and there are constants $\ell_1 = \ell_1(p)$,
$\ell_2=\ell_2(p)
> 0$ (see (\pp{3.102}) below), and a constant $c=c(p)$, such that

$$\|(1-C_{v_\theta})^{-1}\|_{L^p\to L^p} \le \frac{c(1+\lambda)^{\ell_1}}{(1-\rho)^{\ell_2}}\tag{\ppp{1.49}}$$
uniformly for all $x,t\in \bb R$.

\endproclaim
\demo{Remark \ppp{1.50}} Using \pp{1.46}(b), one easily verifies
that estimate \pp{1.49} remains true for $1<p<2$.
\enddemo

Theorem \pp{1.48} plays a crucial role, in particular, in
analyzing the long-time behavior of solutions of the perturbed NLS
equation $iq_t + q_{xx} - 2|q|^2q - \vp|q|^sq=0$, $s>2$ and
$\vp>0$ (see [DZ3], [DZ4], [DZ5]). But the result is also of
independent interest. In the linear case, if $b\in L^\infty(\bb
R)$ and $\phi(z)$ is real valued, then by (\pp{1.2})

$$\|C^-(fbe^{i\phi})\|_{L^p} \le c\|b\|_{L^\infty} \|f\|_{L^p},\tag\ppp{1.51}$$
where the bound is independent of $\phi$. Theorem (\pp{1.48}), in
which the bound $c$ is independent of the multiplier
$e^{i\theta}$, should be viewed as a non-linear version of such
estimates. Indeed if $\rho$ is sufficiently small then we can
expand $(1-C_{v_\theta})^{-1}$ in a Neumann series and a bound of
the form (\pp{1.49}), independent of $e^{i\theta}$, follows from
(\pp{1.51}).

In view of Proposition \pp{1.12}, Theorem \pp{1.48} is equivalent
to the following result.

\proclaim{Theorem \ppp{1.52}}

Suppose $r\in H^1_1, \|r\|_{H^1}\le \lambda, \|r\|_{L^\infty} \le
\rho <1$. Let $\Sigma=\bb R$ and $v=v_\theta$ as in (\pp{1.47}).
Then for $x,t\in \bb R$ and for any $2\le p <\infty$, IRHP$_{L^p}$
has a unique solution $M_\pm$ for any $F$ and there exist $\ell_1
= \ell_1(p)$, $\ell_2  = \ell_2(p)>0$, and a constant $c=c(p)$,
such that

$$\|M_\pm\|_{L^p} \le  \frac{c(1+\lambda)^{\ell_1}}{(1-\rho)^{\ell_2}} \|F\|_{L^p}\tag\ppp{1.53}$$
uniformly for all $x,t\in \bb R$.

\endproclaim

This equivalence allows us to prove (\pp{1.49}) by taking advantage
of the analyticity properties of $M_\pm$. In particular, we are
able to deform the RHP around the stationary phase point $z_0 =
x/2t$ of $e^{i\theta}$, in the spirit of the non-linear
steepest---descent method introduced by the authors in [DZ1]. As
we will see, this deformation plays a crucial role in the
analysis.

\noindent In particular, for the factorization $\left(\matrix
1&-re^{i\theta}\\ 0&1\endmatrix\right)^{-1} \left(\matrix 1&0\\
-\bar re^{-i\theta}&1\endmatrix\right)$ of $v_\theta$,

$$w_\theta = (w^-_\theta,w^+_\theta)  = \left(\left( \matrix 0&re^{i\theta}\\ 0&0\endmatrix\right),\hfil\break \left(\matrix 0&0\\ -\bar re^{-i\theta}&0\endmatrix\right)\right),$$
and it follows that $(1-C_{w_\theta})^{-1}$ obeys the same bound
as $(1-C_{v_\theta})^{-1}$ (cf.\ (\pp{1.49})),

$$\|(1-C_{w_\theta})^{-1}\|_{L^p} \le  \frac{c(1+\lambda)^{\ell_1}}{(1-\rho)^{\ell_2}}\tag\ppp{1.54}$$
uniformly for all $x,t\in \bb R$.

\remark{Remark \ppp{1.55}}

As we will see in Section 3, a simple argument shows that
$(1-C_{v_\theta})^{-1}$ exists as a bounded operator in $L^2(\bb
R)$ and satisfies the bound $\|(1-C_{v_\theta})^{-1}\|_{L^2} \le
c/(1-\rho)$, uniformly for $x,t\in\bb R$. As $v_\theta(z)$ is
continuous in $z$, it follows by general Fredholm arguments (cf.\
[CG]) that for any $x$ and $t$, $(1-C_{v_\theta})^{-1}$ exists as
a bounded operator in $L^p(\bb R)$ for  any $1<p<\infty$, but the
bound may depend on $x$ and $t$. The point here is that, by the
Conjugation Lemma \pp{1.39},  for any $1<p<\infty$, the bound may
be chosen uniformly in $x$ and $t$, as in (\pp{1.49}) and Remark
\pp{1.50}.

The paper is organized as follows. In Section 2 we prove the
Conjugation Lemma \pp{1.39}). In Section~3 we use Lemma \pp{1.39}
to prove Theorem~\pp{1.48}. As noted above, steepest descent
methods play a crucial role (see [DZ1], [DIZ], [DZ2]).

Let $M_\pm\in \partial C(L^p)$ solve the IRHP$_{L^p(\bb R)}$, $M_+ = M_-v_\theta + F$.
If $M_\pm = C^\pm h$, set $\widetilde M_\pm = C^\pm \tilde h$, where $\tilde h(z) = \ovl{h(-z)}$.
(Thus $M(z)  = \ovl{\widetilde M(-\bar z)}$ for the extensions of $M_\pm$,
$\widetilde M_\pm$ off $\bb R$, respectively.) Then a simple computation shows
that $\widetilde M_\pm$ solves the IRHP$_{L^p(\bb R)}$ $\widetilde M_+ =
\widetilde M_- \tilde v_\theta + \widetilde F$, where

$$\tilde v_\theta(z) = \left(\matrix
1-|\tilde r(z)|^2&\tilde r(z) e^{-i\theta(-z)}\\
-\ovl{\tilde r(z)} e^{i\theta(-z)}&1\endmatrix\right)$$ and
$\tilde r(z) = \ovl{r(-z)}, \widetilde F(z) = \ovl{F(-z)}$. As
$\tilde r\in H^1_1$, $\|\tilde r\|_{H^1} = \|r\|_{H^1}$, $\|\tilde
r\|_{L^\infty} = \|r\|_{L^\infty}$, $\|\tilde F\|_{L^p} =
\|F\|_{L^p}$, and as $-\theta(-z) = xz  +tz^2$, it follows from
Proposition~\pp{1.12} that we only need to prove Theorem~\pp{1.48}
for $t\ge 0$. In the text that follows we will assume that $t\ge
0$ without further comment.

\endremark

\remark{Remark on Notation \ppp{1.56}}

If $A = (a_{ij})$ is an $\ell\times m$ matrix, it is convenient in
the remainder of the paper to fix the matrix norm, $|A| \equiv
(\Sigma_{i,j}|a_{ij}|^2)^{\frac12} = (\text{tr } A^*A)^{\frac12}$.
We say $A(z) = (a_{ij}(z))$ is in $L^p(\Sigma)$ for some contour
$\Sigma\subset \bb C$ if each of the entries $a_{ij}(z) \in
L^p(\Sigma)$ and we define $\|A\|_{L^p(\Sigma)} \equiv
\|~|A|~\|_{L^p(\Sigma)}$.

\endremark


\head 2. Proof of the Conjugation Lemma\endhead The proof is in
steps.


\bigskip

\n {\bf Step 1.} As $(1-C_{\breve v})^{-1}$ exists in $L^2$ and
$L^p$, and as $\breve v-I\in L^2\cap L^p$, it follows from
(\pp{1.25})(\pp{1.14}) that $\breve \phi_\pm$, the solution of the
normalized RHP $(\Gamma,\breve v)$, exists in $I+L^p\cap L^2$ and
$\breve\phi_-$ is given by

$$\breve\phi_- = I+(1-C_{\breve v})^{-1} C^-(\breve v-I).\tag\ppp{2.1}$$
Now consider the IRHP$_{L^p}$ on $\Gamma$

$$\breve M_+ = \breve M_- \breve v + F,\quad \breve M_\pm \in \partial C(L^p),\qquad F\in L^p.$$
Then

$$\breve M_+ \breve \phi^{-1}_+ = \breve M_-\breve\phi^{-1}_- + F\breve \phi^{-1}_+$$
and by the Plemelj formula,

$$\breve M_- \breve \phi^{-1}_- = C^-(F\breve \phi^{-1}_+).$$
Thus we obtain by (\pp{1.14}),

$$(C^-(F\breve\phi^{-1}_+))\breve\phi_- = \breve M_- = (1-C_{\breve v})^{-1} C^- F,$$
which implies that $(C^-(\cdot\ \breve\phi^{-1}_+))\breve \phi_-$
is bounded from $L^p\to L^p$ and

$$\|(C^-(\cdot\ \breve\phi^{-1}_+))\breve\phi_-\|_{L^p (\Gamma)} \le c\|(1-C_{\breve v})^{-1}\|_{L^p(\Gamma)
\to L^p(\Gamma)}.\tag\ppp{2.2}$$

{\bf Step 2.} As $(1-C_v)^{-1}$ exists in $L^2$, and as $v-I\in L^2$, it follows as above that $\phi_\pm$,
the solution of the normalized RHP $(\Gamma,v)_{L^2}$ exists and $\phi_-$ is given by

$$\phi_- = I+(1-C_v)^{-1} C^-(v-I).\tag\ppp{2.3}$$
Now let $M_\pm$ solve the IRHP$_{L^2}$ on $\Gamma$

$$M_+ = M_-v + G,\quad M_\pm \in \partial C(L^2)\tag\ppp{2.4}$$
with $G\in L^p\cap L^2\subset L^2$. Then arguing as above,

$$M(z) = (C(G\phi^{-1}_+))(z)\phi(z),\qquad z\in \bb C\backslash \Gamma,$$
where $M(z),\phi(z)$ are the extensions of $M_\pm, \phi_\pm$ off
$\Gamma$ respectively. Write

$$G\phi^{-1}_+ = G + G(\phi^{-1}_+-I) \in L^p(\Gamma) + L^q(\Gamma)$$
where $\frac1q = \frac12 + \frac1p$ i.e. $1<q = \frac{2p}{2+p} < 2
< p$. But by standard estimates $C_{\Gamma\to\Gamma'}$ is bounded
from $L^{p'}(\Gamma) \to L^{p'} \cap L^\infty(\Gamma')$ for any
$p'>1$, and hence

$$\align
\|M\|_{L^p(\Gamma')} &\le c\|\phi\|_{L^\infty(\Gamma')} (\|G\|_{L^p(\Gamma)} +
\|G\|_{L^p(\Gamma)} \|\phi^{-1}_+ - I\|_{L^2(\Gamma)})\tag\ppp{2.5}\\
&\le c(1+ \|\phi_+ - I\|_{L^2(\Gamma)} + \|\phi_- - I\|_{L^2(\Gamma)})(1 + \|\phi^{-1}_+ - I\|_{L^2(\Gamma)}) \|G\|_{L^2(\Gamma)}\\
&\le c\|(1-C_v)^{-1}\|_{L^2(\Gamma)\to L^2(\Gamma)}^2\mnorm
v^2(1+\|v-I\|_{L^2(\Gamma)}^2)\|G\|_{L^p(\Gamma)}.
\endalign$$
Here we used the fact that $\text{dist}(\Gamma,\Gamma')>0$ to
estimate $\|\phi\|_{L^\infty(\Gamma')}$ in terms of
$\|\phi_\pm-I\|_{L^2(\Gamma)}$.

{\bf Step 3.} Again consider the IRHP$_{L^2}$ (\pp{2.4}).
Inserting the relations $v = R^{-1}_- \breve vR_+$ and $\breve v =
(\breve \phi_-)^{-1} \breve \phi_+$, we obtain

$$M_+R^{-1}_+ \breve\phi^{-1}_+ = M_-R^{-1}_- \breve\phi^{-1}_- + GR^{-1}_+
\breve\phi^{-1}_+ \text{ on } \Gamma.\tag\ppp{2.6}$$ Once again
by the Plemelj formula, on $\Gamma$

$$\align
(MR^{-1}\breve\phi^{-1})_- &= C^-_{\Gamma_{\text{ext}} \to
\Gamma}((MR^{-1} \breve\phi^{-1})_+  - (MR^{-1} \breve\phi^{-1})_-)\\
&= C_{\Gamma'\to \Gamma}(MR^{-1}_+ \breve\phi^{-1} - MR^{-1}_-
\breve\phi^{-1}) + C^-_{\Gamma\to\Gamma}(GR^{-1}_+
\breve\phi^{-1}_+).
\endalign$$
Here we have used (\pp{2.6}) and the fact that $M$ and
$\breve\phi$ are analytic across $\Gamma'$. Thus

$$M_- = [C_{\Gamma'\to\Gamma}(MR^{-1}_+ \breve\phi^{-1} -
MR^{-1}_-\breve\phi^{-1})] \breve\phi_- R_- +
(C^-_{\Gamma\to\Gamma}(GR^{-1}_+ \breve\phi^{-1}_+))
\breve\phi_-R_- = \text{I}  + \text{II}.\tag\ppp{2.7}$$ Now as $
\text{dist}(\Gamma,\Gamma') > 0$,

$$
\|C_{\Gamma'\to\Gamma} (MR^{-1}_+ \breve\phi^{-1} - MR^{-1}_-
\breve\phi^{-1})\|_{L^p\cap L^\infty(\Gamma)} \le
c\|R^{-1}\|_{L^\infty(\bb C\backslash\Gamma_{\text{ext}})}
\|M\|_{L^p(\Gamma')} \|\breve\phi\|_{L^\infty(\Gamma')}.
\tag\ppp{2.8}$$ As in Step 2,

$$\|\breve\phi\|_{L^\infty(\Gamma')} \le c\mnorm {\breve v}\|_{L^\infty(\Gamma)}
\|(1-C_{\breve v})^{-1}\|_{L^2(\Gamma)\to L^2(\Gamma)} (1 +
\|\breve v-I\|_{L^2(\Gamma)}),\tag\ppp{2.9}$$ and hence, using
(\pp{2.5}), we obtain

$$\align
&\|C_{\Gamma'\to\Gamma} (MR^{-1}_+ \breve\phi^{-1} - MR^{-1}_-
\breve\phi^{-1})\|_{L^p\cap L^\infty(\Gamma)} \tag\ppp{2.10}\\
&\quad \le c\|R^{-1}\|_{L^\infty(\bb C\backslash \Gamma_{\text{ext}})} \mnorm{\breve v}\mnorm v^2 \|(1-C_{\breve v})^{-1}\|_{L^2(\Gamma)\to L^2(\Gamma)} \|(1-C_v)^{-1}\|^2_{L^2(\Gamma)\to L^2(\Gamma)}\\
&\qquad \times (1+\|\breve v-I\|_{L^2(\Gamma)})
(1+\|v-I\|_{L^2(\Gamma)})^2 \|G\|_{L^p(\Gamma)}.
\endalign$$
Then

$$\|\roman I\|_{L^p(\Gamma)} \le \|C_{\Gamma'\to\Gamma} (MR^{-1}_+
\breve\phi^{-1} - MR^{-1}_- \breve\phi^{-1})\|_{L^p\cap
L^\infty(\Gamma)} \|R\|_{L^\infty}(1 + \|(1-C_{\breve v})^{-1}
\|_{L^p\to L^p} \|\breve v-I\|_{L^p})\tag\ppp{2.11}$$ and, by
(\pp{2.2}),

$$\|\text{II}\|_{L^p(\Gamma)} \le c\|(1-C_{\breve v})^{-1} \|_{L^p(\Gamma)\to L^p(\Gamma)}
 \|R\|_{L^\infty} \|R^{-1}\|_{L^\infty} \|G\|_{L^p}\tag\ppp{2.12}$$
which yields an $L^p(\Gamma)$ bound for $M_-$.

Finally, from (\pp{1.13}), $(1-C_v)^{-1}F = M_- + F$, where
$M_\pm$ solve IRHP$_{L^2}$ with $G = F(v-I)$. For $F\in L^p\cap
L^2$, and hence $G\in L^p\cap L^2$, we then obtain from
(\pp{2.11})(\pp{2.12}), together with (\pp{1.23},\pp{1.24}) and
their analogs for $\breve v$, and also the interpolation estimate
$\|v-I\|_{L^p}\le c(\|v-I\|_{L^\infty}+\|v-I\|_{L^2})\le c\mnorm
v(1+\|v-I\|_{L^2})$, a bound of the form $\|(1-C_v)
^{-1}F\|_{L^p(\Gamma)} \le c^\#_v\|F\|_{L^p(\Gamma)}$ with
$c^\#_v$ as in (\pp{1.45}). The result then follows by
density.$\qquad \square$


The constant $c^\#_v$ in (\pp{1.44}) depends on $\Gamma$ and
$\Gamma'$ and the distance $\text{dist}(\Gamma,\Gamma')$ between
them. We are interested in particular in applying the Conjugation
Lemma~\pp{1.39} for contours of the type that appear in
Figure~\pp{1.38}. For such contours it is easy to show that for
$1<q\le p\le \infty$,

$$\|C^\pm\|_{L^q(\Sigma')\to L^p(\Sigma)},\quad \|C^\pm\|_{L^q(\Sigma)\to L^p(\Sigma')}
\le \frac{c_q}{(\text{dist }
\Gamma,\Gamma')^{\frac1q-\frac1p}}.\tag\ppp{2.13}$$ We are
particularly interested  in the case where
$\text{dist}(\Gamma,\Gamma')$ is bounded, say
$\text{dist}(\Gamma,\Gamma') \le 1$. Keeping track of the
constants, and using $\text{dist}(\Gamma,\Gamma')\le 1$, we see
that the constant $c$ in (\pp{2.5}) should be replaced by
$c/\text{dist}(\Gamma,\Gamma')$, and that $c$ in (\pp{2.10})
should be replaced by $c/(\text{dist}(\Gamma,\Gamma'))^{3/2+1/p}$.
It follows that in (\pp{1.44}) we should replace

$$c^\#_v \to c^\#_v/(\text{dist}(\Gamma,\Gamma'))^{3/2+1/p}\tag\ppp{2.14} $$
in the case $\text{dist}(\Gamma,\Gamma')\le 1.$


\head 3. Proof of Theorem \pp{1.48}\endhead

\bigskip

\demo{Notational Remark} In this section all jump matrices $v$ are
$2\times2$ with determinant 1. Thus if $v=\pmtwo abcd$, then
$v^{-1}=\pmtwo d{-b}{-c}a$. It follows that we can replace $\mnorm
v$ by $\|v\|_{L^\infty}$ in all the estimates in Section 2, and in
particular, in (\pp{1.45}). We will make this replacement
systematically in this section without further comment.
\enddemo

As indicated in Remark \pp{1.55}, the proof of Theorem~\pp{1.48} in
the case $p=2$ follows by a simple argument. Indeed for the
operator $C_{w_\theta}$ in (\pp{1.54}), where

$$w_\theta = (w^-_\theta, w^+_\theta) = \left(\left(
\matrix 0&re^{i\theta}\\ 0&0\endmatrix\right), \left(\matrix 0&0\\
-\bar re^{-i\theta}&0\endmatrix\right)\right),$$ we have for $h =
(h_{ij})_{1\le i,j\le 2}$,
$$C_{w_\theta}h = \left(\matrix
C^-(-h_{12}\bar re^{-i\theta})&C^+(h_{11} re^{i\theta})\\
C^-(-h_{22}\bar re^{-i\theta})&C^+(h_{21}re^{i\theta})
\endmatrix\right).$$ But under Fourier transform $C^+$ (resp.\
$-C^-$) is just multiplication by the characteristic function of
$(0,\infty)$ (resp. $(-\infty,0)$) and hence $\|C^\pm\|_{L^2(\bb
R)} = 1$. Thus (recall Remark on Notation \pp{1.56}),
$\|C_{w_\theta}h\|_{L^2} \le \|r\|_{L^\infty} \|h\|_{L^2}$ and
hence $\|C_{w_\theta}\|_{L^2} \le \|r\|_{L^\infty}$. Thus for
$\|r\|_\infty \le \rho < 1$, it follows that
$(1-C_{w_\theta})^{-1}$, and hence $(1-C_{v_\theta})^{-1}$ (use
(\pp{1.16})), exist and are  uniformly bounded in $L^2(\bb R)$ for
all $x,t\in \bb R$,

$$\|(1-C_{w_\theta})^{-1}\|_{L^2(\bb R)\to\rb},\  \|(1-C_{v_\theta})^{-1}\|_{L^2(\bb R)\to L^2(\rb)}
\le \frac c{1-\rho}.\tag\ppp{3.1}$$
Note that the bound depends only on $\rho$ and not on the $H^1$
norm of $r$.

We now give a second proof of (\pp{3.1}) by a more general method
that will be useful at various points in the calculations that
follow.

\proclaim{Proposition \ppp{3.2}}

Suppose $r\in L^\infty(\bb R)$ and $\|r\|_\infty\le \rho <1$. Let
$v = \left(\matrix 1-|r|^2&r\\ -\bar r&1\endmatrix\right)$, and
let $C_v$ be the associated Cauchy operator as in (\pp{1.5}). Then
$(1-C_v)^{-1}$ exists in $L^2(\bb R)$ and

$$\|(1-C_v)^{-1}\|_{L^2(\bb R)\to L^2(\bb R)} \le \frac{c}{1-\rho}.\tag\ppp{3.3}$$

\endproclaim

\demo{Proof}

Assume first that $(1-C_v)^{-1}$ exists. Then by (1.11), $(1-C_v)^{-1} f = (M_-+f)$ for $f\in L^2$, where $M$ solves the IRHP$_{L^2}$, $M_+ = M_-v + f(v-1)$, $M_\pm \in \partial C(L^2)$. Now by a simple contour argument

$$\int_{\bb R} M_-(z) M^*_+(z)dz = 0,$$

and hence

$$\int_{\bb R} M_-(z) v^*(z) M^*_-(z)dz = \int_{\bb R} M_-(z) (I-v^*(z)) f^*(z)dz.$$

Taking conjugates and adding we obtain

$$\int_{\bb R} M_-(z)(v(z) + v^*(z)) M^*_-(z)dz = \int_{\bb R} [M_-(z)(I-v^*(z))f^*(z) + f(z)(I-v(z))M^*_-(z)]dz.\tag\ppp{3.4}$$
But $v(z) + v(z)^* = 2 \left(\matrix 1-|r(z)|^2&0\\
0&1\endmatrix\right)$, and it then follows directly from
(\pp{3.4}) that $\|M_-\|_{L^2} \le \frac{c}{1-\rho} \|f\|_{L^2}$,
which implies in turn the bound (\pp{3.3}) for $(1-C_v)^{-1}$. Set
$r_\gamma = \gamma r$, $0\le \gamma \le 1$. We have $v_\gamma-I=
\left(\matrix -\gamma^2|r|^2&\gamma r\\ -\gamma\bar
r&0\endmatrix\right)$ and hence $1-C_{v_\gamma} =
1-C(\cdot(v_\gamma-I))$ is invertible for small $\gamma$. However
$\|r_\gamma\|_{L^\infty} \le \rho$ for all $0 \le \gamma\le 1$,
and hence $(1-C_{v_\gamma})^{-1}$ must satisfy the bound
(\pp{3.3}) whenever it exists. By an elementary continuity
argument it then follows that $(1-C_{v_\gamma})^{-1}$ exists and
satisfies (\pp{3.3}) for all $0\le \gamma\le 1$.$\qquad \square$

\enddemo

In order to control $(1-C_{v_\theta})^{-1}$ in $L^p$, $2<p<\infty$, uniformly for $x,t\in\bb R$,
we must control the solution $M_\pm$ of the IRHP$_{L^p}$ with jump matrix $v_\theta$. Following
the steepest descent method introduced in [DZ1], and applied to the NLS equation in [DIZ], [DZ2],
we expect the IRHP$_{L^p}$ to ``localize'' near the stationary phase point $z_0 = x/2t$ for
 $\theta = xz-tz^2$, $\theta'(z_0) = 0$. Furthermore, the signature table for $\text{Re } i\theta$

\vskip0.5cm

\font\thinlinefont=cmr5

\centerline{\beginpicture \setcoordinatesystem units
<.7000cm,.7000cm>
\setshadesymbol ({\thinlinefont .}) \setlinear
%
%
\linethickness= 0.500pt \setplotsymbol ({\thinlinefont .})
\putrule from  2.540 21.590 to  7.620 21.590
%
%
\linethickness= 0.500pt \setplotsymbol ({\thinlinefont .})
\putrule from  5.080 24.130 to  5.080 19.050
%
%
\put{$z_0$} [lB] at  5.239 21.114
%
%
\put{$\re i\theta<0$} [lB] at  2.651 22.860
%
%
\put{$\re i\theta>0$} [lB] at  5.7 22.860
%
%
\put{$\re i\theta>0$} [lB] at  2.651 20.161
%
%
\put{$\re i\theta<0$} [lB] at  5.7 20.161 \linethickness=0pt
\putrectangle corners at  2.515 24.155 and  7.645 19.025
\endpicture}

\centerline{Figure \ppp{3.5}. Signature table for $\text{Re }
i\theta$}

\bigskip

\n should play a crucial role. The basic idea of the method is to
deform the contour $\Gamma=\bb R$ so that the exponential factors
$e^{i\theta}$ and $e^{-i\theta}$ are exponentially decreasing, as
dictated by Figure~\pp{3.5}.
 In order to make these deformations we must separate the factors $e^{i\theta}$ and $e^{-i\theta}$
 algebraically, and this is done using the upper/lower and lower/upper factorizations of $v_\theta$,

$$v_\theta = \left(\matrix 1&re^{i\theta}\\ 0&1\endmatrix\right)
\left(\matrix 1&0\\ -\bar re^{-i\theta}&1\endmatrix\right) =
\left(\matrix 1&0\\ \frac{-\bar
r}{1-|r|^2}e^{-i\theta}&1\endmatrix\right) \left(\matrix
1-|r|^2&0\\ 0&\frac1{1-|r|^2}\endmatrix\right) \left(\matrix 1&
\frac{r}{1-|r|^2} e^{i\theta}\\
0&1\endmatrix\right).\tag\ppp{3.6}$$

The upper/lower factorization is appropriate for $z>z_0$ and the lower/upper factorization is appropriate for $z<z_0$.
 The diagonal terms in the lower/upper factorization can be removed by conjugating $v_\theta$,

$$\breve v_\theta = \delta^{\sigma_3}_- v_\theta \delta^{-\sigma_3}_+,\quad \sigma_3 =
\left(\matrix 1&0\\ 0&-1\endmatrix\right) = \text{ third Pauli
matrix,}\tag\ppp{3.7}$$ by the solution $\delta_\pm$ of the
scalar, normalized RHP $(\bb R_- + z_0,1-|r|^2)_{L^2}$,

$$\left\{\matrix \delta_+ = \delta_-(1-|r|^2),\hfill& z\in \bb R_- + z_0,\hfill\\
\delta_\pm - 1\in \partial C(L^2),\hfill\endmatrix\right.
\tag\ppp{3.8}$$ where the contour $\bb R_- + z_0$ is oriented from
$-\infty$ to $z_0$. The properties of $\delta$ can be read off
from the following elementary proposition, whose proof is left to
the reader.

\proclaim{Proposition \ppp{3.9}}

Suppose $r\in L^\infty(\bb R)\cap L^2(\bb R)$ and
$\|r\|_{L^\infty} \le \rho <1$. Then the solution $\delta_\pm$ of
the scalar normalized RHP (\pp{3.8}) exists and is unique and is
given by the formula

$$\delta_\pm(z) = e^{C^\pm_{\bb R_- + z_0} \log 1-|r|^2} = e^{\frac1{2\pi i}
\int^{x_0}_{-\infty} \frac{\log(1-|r(s)|^2)}{s-z_\pm}},\qquad z\in
\bb R.\tag\ppp{3.10}$$ The extension $\delta$ of $\delta_\pm$ off
$\bb R_- +z_0$ is given by

$$\delta(z) = e^{C_{\bb R_- + z_0} \log(1-|r|^2)} = e^{\frac1{2\pi i}
\int^{z_0}_{-\infty} \frac{\log(1-|r(s)|^2)}{s-z}ds},\qquad
z\in\bb C\backslash(\bb R_- + z_0),\tag\ppp{3.11}$$ and satisfies
for $z\in \bb C \backslash(\bb R_- + z_0)$,

$$\gather
\delta(z)\ovl{\delta(\bar z)} = 1,\tag\ppp{3.12}\\
(1-\rho)^{\frac12} \le (1-\rho^2)^{\frac12} \le |\delta(z)|,
|\delta^{-1}(z)| \le (1-\rho^2)^{-\frac12} \le
(1-\rho)^{-\frac12},\tag\ppp{3.13}
\endgather$$
and

$$|\delta^{\pm 1}(z)| \le 1\quad \text{for}\quad\pm \text{\rm Im } z>0.\tag\ppp{3.14}$$
For real $z$,

$$\gather
|\delta_+(z)\delta_-(z)| = 1 \quad \text{ and, in particular,}\quad |\delta(z)| =
1\quad \text{for}\quad z>z_0,\tag\ppp{3.15}\\
|\delta_+(z)| = |\delta^{-1}_-(z)| = (1-|r(z)|^2)^{\frac12},
\qquad z<z_0,\tag\ppp{3.16}
\endgather$$
and

$$\Delta \equiv \delta_+\delta_- = e^{\frac1{i\pi} \text{ P.V. } \int^{z_0}_{-\infty}
\frac{\log(1-|r(s)|^2)}{s-z}ds},\quad \text{where P.V.\ denotes
the principal value.}\tag\ppp{3.17}$$ Also $|\Delta| =
|\delta_+\delta_-| = 1$

$$\|\delta_\pm - 1\|_{L^2(dz)} \le \frac{c\|r\|_{L^2}}{1-\rho}.\tag\ppp{3.18}$$

\endproclaim

We obtain the following factorization for $\breve v_\theta$

$$\align
\breve v_\theta = \breve v^{-1}_{\theta-} \breve v_{\theta+} &=
\left(\matrix 1&re^{i\theta} \delta^2\\ 0&1\endmatrix\right)
\left(\matrix 1&0\\ -\bar re^{-i\theta} \delta^{-2}&1\endmatrix\right),\qquad z>z_0, \tag\ppp{3.19}\\
\breve v_\theta = \breve v^{-1}_{\theta-} \breve v_{\theta+} &=
\left(\matrix 1&0\\ \frac{-\bar re^{-i\theta}
\delta^{-2}_-}{1-|r|^2}&1\endmatrix\right)
\left(\matrix1&\frac{re^{i\theta}\delta^2_+}{1-|r|^2}\\
0&1\endmatrix\right),\qquad z<z_0.\tag\ppp{3.20}
\endalign$$
Using Figure \pp{3.5} we observe the crucial fact that the
analytic continuations to $\bb C_+$ of the exponentials in the
factors on the right in (\pp{3.19}) and (\pp{3.20}), are
exponentially decreasing, and the same is true for the
exponentials on the left, when continued to $\bb C_-$.

For later reference, observe that (\pp{3.20}) can also be written
in the form

$$\breve v_\theta = \left(\matrix 1&0\\ -\bar re^{-i\theta} \delta^{-1}_+
\delta^{-1}_-&1\endmatrix\right) \left(\matrix
1&re^{i\theta}\delta_+\delta_-\\ 0&1\endmatrix\right), \qquad
z<z_0.\tag\ppp{3.21}$$

Now clearly

$$\gather
\breve M_\pm \text{ solves IRHP$_{L^2}(\bb R, \breve v)$ with inhomogeneous term } \breve F\\
\Leftrightarrow\\
M_\pm = \breve M_\pm \delta^{\sigma_3}_\pm \text{ solves
IRHP$_{L^p}$ with inhomogeneous term } F = \breve F
\delta^{\sigma_3}_+
\endgather$$
and so to control $(1-C_{v_\theta})^{-1}$, it is sufficient to
control $(1-C_{\breve v_\theta})^{-1}$.

The $L^p$ bound on $(1-C_{\breve v_\theta})^{-1}$, and hence on
$(1-C_{v_\theta})^{-1}$, in the general case will be inferred
eventually from the following model problem, with the aid of the
Conjugation Lemma~\pp{1.39}.

\medskip

\n (\ppp{3.22})\ {\bf Model Problem:}

\itemitem{(i)} Suppose $x=0$ so that $z_0=0$ and $\theta =-tz^2$

\itemitem{(ii)} $r(z) = \frac{r(0)}{1+iz}$, where $|r(0)| \le \rho <1$.

\proclaim{Proposition \ppp{3.23}}

Let $0<\beta<1/2$. Then for the model problem (\pp{3.22}),

$$\|(1-C_{\breve v_\theta})^{-1}\|_{L^p\to L^p} \le \frac{c}{(1-\rho)^{2+5\beta}},\qquad 2<p<\infty,\tag\ppp{3.24}$$
where $c$ is independent of $t>0$ and $\rho<1$, and depends only
on $\beta$ on $p$.

\endproclaim

We will prove this Proposition in a series of steps. Let
$\delta_\pm$ denote the solution of the scalar normalized RHP
(\pp{3.8}) for the model $r$ in (\pp{3.22}).\medskip

\n {\bf Step 1 (analytic continuations).} Observe that for $r =
\frac{r(0)}{1+iz}$ and $z<0$, $\delta_+\delta_- =
\frac{\delta^2_+}{1-|r|^2}  = \frac{(z^2+1)\delta^2_+}{z^2 +
1-|r(0)|^2}$, and so $h \equiv \frac{z-i\sqrt{1 -
|r(0)|^2}}{z+i\sqrt{1-|r(0)|^2}} \delta_+\delta_-$ has an analytic
continuation from $\bb R_-$ to $C_+$. For $z<0$, $|h(z)| \le 1$ by
(\pp{3.15}), and for $z>0$, $|h(z)| \le \frac1{1-|r(0)|^2} \le
\frac1{1-\rho^2}$, again by (\pp{3.15}). By a standard
Phragm\'en--Lindel\"of argument we conclude that

$$|(\delta_+\delta_-)(z)| \le \left| \frac{z+i\sqrt{1-|r(0)|^2}}{z-i\sqrt{1-|r(0)|^2}}\right|
\frac1{(1-\rho^2)^{1 - \frac{\arg z}{\pi}}}\tag\ppp{3.25}$$ for
$z\in\bb C_+$, $0<\arg z<\pi$. It follows that
$(\delta_+\delta_-)^{-1}(z) = \ovl{(\delta_+\delta_-)(\bar z)}$,
$z<0$, has an analytic continuation from $\bb R_-$ to $\bb C_-$
satisfying

$$|(\delta_+\delta_-)^{-1}(z)| \le \left|\frac{z-i \sqrt{1-|r(0)|^2}}
{z+i\sqrt{1-|r(0)|^2}}\right| \frac1{(1-\rho^2)^{1 + \frac{\arg
z}\pi}}\tag\ppp{3.26}$$ for $-\pi < \arg z < 0$. Also the analytic
continuations of $\delta(z)^{-2}$, $\delta(z)^2$ to $\bb C_+$,
$\bb C_-$ respectively satisfy the bounds

$$|\delta(z)^{-2}|\le \frac1{(1-\rho^2)^{\frac{\arg z}\pi}}, \qquad z\in\bb C_+, \quad 0<\arg z<\pi,\tag\ppp{3.27}$$
and

$$|\delta(z)^2| \le (1-\rho^2)^{\frac{\arg z}\pi},\qquad z\in \bb C_-,\quad -\pi < \arg z < 0.\tag\ppp{3.28}$$
It follows that $\breve v_{\theta+}$ and $\breve v_{\theta-}$ in
(\pp{3.19}) and (\pp{3.21}), have analytic continuations to $\bb
C_+$ and $\bb C_-$ respectively, where they satisfy the bounds
$(0<\beta<\frac12)$

$$\align
|\breve v_{\theta+}(z)-I| &\le \frac{c}{(1-\rho)^\beta}\quad
\text{for } z\in \bb C_+, \arg z = \beta \pi \text{ or } \pi - \beta\pi,\tag\ppp{3.29}\\
|\breve v_{\theta-}(z)-I|&\le \frac{c}{(1-\rho)^\beta} \quad
\text{for } z\in \bb C_-, \arg z = -\beta \pi \text{ or } -\pi +
\beta\pi,\tag\ppp{3.30}
\endalign$$
as indicated in Figure \pp{3.31}

\vskip0.5cm

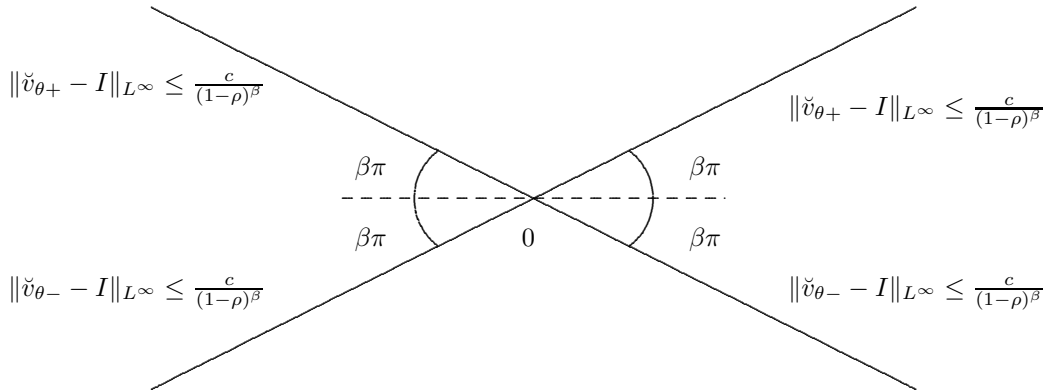
\font\thinlinefont=cmr5 \centerline{\beginpicture
\setcoordinatesystem units <1.00000cm,1.00000cm> \setshadesymbol
({\thinlinefont .}) \setlinear
%
%
\linethickness= 0.500pt \setplotsymbol ({\thinlinefont .})
\circulararc 106.260 degrees from  6.350 22.225 center at  6.826
21.590
%
%
\linethickness= 0.500pt \setplotsymbol ({\thinlinefont .})
\circulararc 106.260 degrees from  8.890 20.955 center at  8.414
21.590
%
%
\linethickness= 0.500pt \setplotsymbol ({\thinlinefont .}) \plot
2.540 24.130 12.700 19.050 /
%
%
\linethickness= 0.500pt \setplotsymbol ({\thinlinefont .}) \plot
2.540 19.050 12.700 24.130 /
%
%
\linethickness= 0.500pt \setplotsymbol ({\thinlinefont .})
\setdashes < 0.1270cm> \plot  5.080 21.590 10.160 21.590 /
%
%
\put{0} [lB] at  7.461 20.955
%
%
\put{$\beta\pi$} [lB] at  5.239 21.907
%
%
\put{$\beta\pi$} [lB] at  9.684 21.907
%
%
\put{$\beta\pi$} [lB] at  5.239 20.955
%
%
\put{$\beta\pi$} [lB] at  9.684 20.955
%
%
\put{$\|\circv_{\theta+}-I\|_{L^\infty}\le\frac
c{(1-\rho)^\beta}$} [lB] at  0.635 23.019
%
%
\put{$\|\circv_{\theta+}-I\|_{L^\infty}\le\frac
c{(1-\rho)^\beta}$} [lB] at 11 22.701
%
%
\put{$\|\circv_{\theta-}-I\|_{L^\infty}\le\frac
c{(1-\rho)^\beta}$} [lB] at  0.635 20.3
%
%
\put{$\|\circv_{\theta-}-I\|_{L^\infty}\le\frac
c{(1-\rho)^\beta}$} [lB] at 11 20.3 \linethickness=0pt
\putrectangle corners at  0.635 24.155 and 12.725 19.025
\endpicture}

\centerline{Figure \ppp{3.31}. Bounds for the continuations of
$\circv_{\theta\pm}$.}\bigskip

\n The constants $c$ in (\pp{3.30}) are independent of $\rho$ and
$t>0$.\medskip

\n {\bf Step 2 (scaling and augmentation).} It is convenient to
scale the RHP as follows:
$$\circv_\theta \to \circv_t(z) \equiv ~\circv_\theta
(z/\sqrt t) = e^{-iz^2 \ad \sigma}\circv(z/\sqrt t).
\tag\ppp{3.32}$$ If $S_t$ denotes the scaling operator $S_tf(z) =
t^{-\frac1{2p}} f(z/\sqrt t)$, then $S_t$ is an isometry from
$L^p({\Bbb R})$ onto $L^p({\Bbb R})$ and
$$\frac1{1-C_{v_\theta}} = S^{-1}_t \left(\frac1{1-
C_{\circv_t}}\right)S_t\tag\ppp{3.33}$$ and hence to prove
(\pp{3.24}) it is enough to replace $\circv_\theta$ with
$\circv_t$. We denote the associated factors of $\circv_t$ by
$\circv_{t\pm}$ and as the bounds in (\pp{3.29})(\pp{3.30}) are
unaffected by scaling, they remain true for $\circv_{t\pm}$ for
all $t>0$.

Consider the IRHP$_{L^p}$ on $\Gamma = {\Bbb R}$,
$$M_+ = M_-\circv_t + F,\qquad F\in L^p({\Bbb R}).
\tag\ppp{3.34}$$ Extend (\pp{3.34}) trivially to the augmented
contour $\Gamma$ in Figure~\pp{3.35} with an opening angle
$\beta\pi$ as in Figure \pp{3.31}

\bigskip\bigskip
\font\thinlinefont=cmr5

\centerline{\beginpicture \setcoordinatesystem units <.7cm,.7cm>
\setshadesymbol ({\thinlinefont .}) \setlinear
%
%
\linethickness= 0.500pt \setplotsymbol ({\thinlinefont .}) \plot
2.540 24.130 12.700 19.050 /
%
%
\linethickness= 0.500pt \setplotsymbol ({\thinlinefont .}) \plot
2.540 19.050 12.700 24.130 /
%
%
\linethickness= 0.500pt \setplotsymbol ({\thinlinefont .})
\putrule from  2.540 21.590 to 12.700 21.590
%
%
\linethickness= 0.500pt \setplotsymbol ({\thinlinefont .})
\putrule from  3.810 21.590 to  4.128 21.590
%
%
\plot  3.873 21.526  4.128 21.590  3.873 21.654 /
%
%
%
\linethickness= 0.500pt \setplotsymbol ({\thinlinefont .})
\putrule from 10.160 21.590 to 10.478 21.590
%
%
\plot 10.224 21.526 10.478 21.590 10.224 21.654 /
%
%
%
\linethickness= 0.500pt \setplotsymbol ({\thinlinefont .}) \plot
3.810 19.685  4.128 19.844 /
%
%
\plot  3.929 19.673  4.128 19.844  3.872 19.787 /
%
%
%
\linethickness= 0.500pt \setplotsymbol ({\thinlinefont .}) \plot
3.493 23.654  3.810 23.495 /
%
%
\plot  3.554 23.552  3.810 23.495  3.611 23.665 /
%
%
%
\linethickness= 0.500pt \setplotsymbol ({\thinlinefont .}) \plot
10.160 22.860 10.478 23.019 /
%
%
\plot 10.279 22.848 10.478 23.019 10.222 22.962 /
%
%
%
\linethickness= 0.500pt \setplotsymbol ({\thinlinefont .}) \plot
9.842 20.479 10.160 20.320 /
%
%
\plot  9.904 20.377 10.160 20.320  9.961 20.490 /
%
%
%
\put{0} [lB] at  7.461 20.955
%
%
\put{$\Omega_0$} [lB] at  7.144 23.336
%
%
\put{$\Omega_0$} [lB] at  6.985 19.367
%
%
\put{$\Omega_+$} [lB] at  3.334 22.384
%
%
\put{$\Omega_+$} [lB] at 10.795 22.384
%
%
\put{$\Omega_-$} [lB] at  3.334 20.479
%
%
\put{$\Omega_-$} [lB] at 10.795 20.637 \linethickness=0pt
\putrectangle corners at  2.515 24.155 and 12.725 19.025
\endpicture}

\centerline {Figure \ppp{3.35}. Augmented contour
$\Gamma$.}\bigskip

\n by setting

$$\breve v_t\equiv I,\quad F\equiv 0 \text{ on } \Gamma/\bb R.\tag\ppp{3.36}$$
Clearly the extension $M$ of $M_\pm$ off $\bb R$ satisfies the
augmented IRHP$_{L^p}$

$$M_+ = M_- \breve v_t + F \text{ on } \Gamma.\tag\ppp{3.37}$$
Set

$$\cases \Phi = I&\text{in $\Omega_0$}\\
\Phi = \breve v_{t\pm}&\text{in
$\Omega_\pm$}\endcases\tag\ppp{3.38}$$ and define

$$\hat v = \Phi_-\breve v_t \Phi^{-1}_+ \text{ on } \Gamma.\tag\ppp{3.39}$$
Note that

$$\hat v = I \text{ on } \bb R.\tag\ppp{3.40}$$

\n {\bf Step 3. (Bound for $\bold{(1-C_{\hat v})^{-1}}$ in ${\bold{ L^2}}(\boldsymbol\Gamma)$).} We prove the following Lemma.

\proclaim{Lemma \ppp{3.41}}

$$\|(1-C_{\hat v})^{-1}\|_{L^2(\Gamma)\to L^2(\Gamma)} \le \frac{c}{(1-\rho)^{1+3\beta}}\tag\ppp{3.42}$$
where $c$ is independent of $\rho<1$ and $t>0$.

\endproclaim

\demo{Proof}

Using the properties of $\delta$ in Proposition \pp{3.9}, one sees
from (\pp{3.19}) and (\pp{3.21}) that $\breve v_t$ has the form
$\left(\matrix 1-|\alpha|^2&\alpha\\
-\bar\alpha&1\endmatrix\right)$ on $\bb R_+$ and $\left(\matrix
1&\alpha\\ -\bar\alpha&1-|\alpha|^2\endmatrix\right)$ on $\bb R_-$
for
 some function $\alpha$ with $\|\alpha\|_{L^\infty(\bb R)} \le \rho$ for all $t>0$, and it follows then from (the proof of)
 Proposition~\pp{3.2} that $(1-C_{\tilde v_t})^{-1}$ exists in $L^2(\bb R)$ and

$$\|(1-C_{\breve v_t})^{-1}\|_{L^2(\bb R)\to L^2(\bb R)} \le \frac{c}{1-\rho},\tag\ppp{3.43}$$
where $c$ is independent of $t>0$. Now consider the equation
$(1-C_{\breve v_t})f =g$ in $L^2$ of the augmented contour
$\Gamma$. As $\breve v_t = I$ on $\Gamma \backslash \bb R$, this
reduces on $\bb R\subset \Gamma$ to the equation $(1-C^-_{\bb
R}(\cdot\ (\breve v_t - I))(f \upharpoonright \bb R) =
(g\upharpoonright \bb R)$ so that $f \upharpoonright \bb R =
(1-C_{\breve v_t})^{-1}_{L^2(\bb R) \to L^2(\bb R)} (g
\upharpoonright  \bb R)$. But then for $z\in \Gamma\backslash \bb
R$, $f(z) = (C_{\bb R}((f\upharpoonright  \bb R)(\breve
v_t-I)))(z) + g(z)$, and hence $\|f\|_{L^2(\Gamma\backslash \bb
R)} \le \frac{c}{1-\rho} \|g\|_{L^2(\Gamma)}$ by (\pp{3.43}),
which thereby extends to $\Gamma$,

$$\|(1-C_{\breve v_t})^{-1}\|_{L^2(\Gamma)\to L^2(\Gamma)} \le \frac{c}{1-\rho},\tag\ppp{3.44}$$
where again $c$ is independent of $t>0$.

Finally consider the IRHP$_{L^2}$ on $\Gamma$,

$$\widehat M_+ = \widehat M_- \hat v + \widehat F(\hat v - I),\qquad \widehat F\in L^2(\Gamma).\tag\ppp{3.45}$$
By (\pp{3.39}), $\widehat M_+\Phi_+ = \widehat M_-\Phi_- \breve
v_t + \widehat F(\hat v-I)\Phi_+$, and so by (\pp{1.14}),

$$\align
\|\widehat M_-\|_{L^2(\Gamma)} &\le \|(1-C_{\breve v_t})^{-1}
C^-(\widehat F(\hat v-I) \Phi_+)\|_{L^2(\Gamma)} \|\Phi^{-1}_-\|_{L^\infty(\Gamma)}\\
&\le c\|(1-C_{\breve v_t})^{-1}\|_{L^2\to L^2} \|\widehat F\|_{L^2} \|\hat v-I\|_{L^\infty} \|\Phi_+\|_{L^\infty} \|\Phi^{-1}_-\|_{L^\infty}\\
&\le \frac{c}{(1-\rho)^{1+3\beta}} \|\widehat F\|_{L^2}
\endalign$$
by (\pp{3.38}), (\pp{3.40}), (\pp{3.29}), (\pp{3.30}), and
(\pp{3.43}). But then (\pp{3.42}) follows from
(\pp{1.13}).$\qquad\square$

\enddemo

\n {\bf Step 4 (Bound for $\bold{(1-C_{\breve v_t})^{-1}}$ in
$\bold{L^p}$).} As $(1-C_{\hat v_t})^{-1}$ exists in
$L^2(\Gamma)$, we know that the IRHP$_{L^2}$ (\pp{3.37}) has a
solution $M_\pm \in \partial C(L^2(\Gamma))$ for $F\in L^p\cap
L^2\subset L^2(\Gamma)$. Inserting (\pp{3.39}) into (\pp{3.37})
for such $F\in L^p\cap L^2(\Gamma)$, and using (\pp{1.14}), we see
that

$$M_-\Phi^{-1}_- = C^-(F\Phi^{-1}_+) + (1-C_{\hat v})^{-1}
C_{\hat v}(C^-(F\Phi^{-1}_+)) \equiv \text{I} +
\text{II}.\tag\ppp{3.46}$$ Now

$$\align
\|\text{II}\|_{L^2(\Gamma)} &\le c\|(1-C_{\hat v})^{-1}\|_{L^2\to L^2}
\|\hat v-I\|_{L^{p'}} \|F\|_{L^p} \|\Phi^{-1}_
+\|_{L^\infty}\tag \ppp{3.47}\\
&\le \frac{c}{(1-\rho)^{1+5\beta}} \|F\|_{L^p}, \ \ 1/p'+1/p=1/2,
\endalign$$
where we have used (\pp{3.38})--(\pp{3.40}), (\pp{3.29}),
(\pp{3.30}), together with (\pp{3.42}). Note that the exponential
factors $e^{\pm iz^2}$ are suitably disposed with respect to the
signature table in Figure~\pp{3.5} and play an essential role in
ensuring that the $L^{p'}$ norm of $\hat v-I$ is bounded uniformly
in $t$. This fact is at the heart of the utility of the steepest
descent method in proving the uniform $L^p$ bounds on
$(1-C_{\breve v_t})^{-1}$, and hence, eventually, the desired
$L^p$ bounds on $(1-C_{v_\theta})^{-1}$.

Write $\text{II} = C_{\hat v}(C^-(F\Phi^{-1}_+)) + C_{\hat v} \text{II}$. Since $\hat v = I$ on $\bb R$,
 we obtain by (\pp{3.47})

$$\|\text{II}\|_{L^p(\bb R)\backslash [-1,1])} \le \frac{c}{(1-\rho)^{2\beta}}
\|F\|_{L^p(\Gamma)} + \frac{c}{(1-\rho)^{1+6\beta}}
\|F\|_{L^p(\Gamma)} \le \frac{c}{(1-\rho)^{1+6\beta}}
\|F\|_{L^p(\Gamma)}.\tag\ppp{3.48}$$ But

$$\|\text{I}\|_{L^p(\bb R)} \le \frac{c}{(1-\rho)^\beta} \|F\|_{L^p},\tag\ppp{3.49}$$
and hence from (\pp{3.46})

$$\|M_-\|_{L^p(\bb R\backslash [-1,1])} \le \frac{c}{(1-\rho)^{1+6\beta}} \|F\|_{L^p(\Gamma)}\tag\ppp{3.50}$$
as $\|\Phi_-\|_{L^\infty(\bb R)} \le c$.

It remains to estimate $\|M_-\|_{L^p(-1,1)}$. Write $\breve v_t$
in the form $\breve v_t = \delta^{\sigma_3}_{t-} \left(\matrix
1&r_t\\ 0&1\endmatrix\right) \left(\matrix 1&0\\ -\bar
r_t&1\endmatrix\right) \delta^{-\sigma_3}_{t+}$, where $r_t(z) =
e^{-iz^2} r(z/\sqrt t)$, $\delta_t(z) = \delta(z/\sqrt t)$.
Inserting this factorization into (\pp{3.34}), we obtain on $\bb
R$

$$M_+ \delta^{\sigma_3}_{t+} \left(\matrix 1&0\\ \bar r_t&1\endmatrix\right) =
M_-\delta^{\sigma_3}_{t-} \left(\matrix 1&r_t\\
0&1\endmatrix\right) + F\delta^{\sigma_3}_{t+} \left(\matrix 1&0\\
\bar r_t&1\endmatrix\right).\tag\ppp{3.51}$$ From the explicit
form of $r$ in (\pp{3.22}) we see that $M_+ \delta^{\sigma_3}_{t+}
\left(\matrix 1&0\\ \bar r_t&1\endmatrix\right)$ has an analytic
continuation to $\bb C_+$,
and $M_-\delta^{\sigma_3}_{t-} \left(\matrix 1&r_t\\
0&1\endmatrix\right)$ has an analytic continuation to $\bb C_-$.
Denote these continuations by $N$ in $\bb C\backslash \bb R$. By
Cauchy's formula for $|z|< 2$, $z\notin \bb R$, we obtain from
(\pp{3.51})

$$N(z) = \oint\limits_{|s|=2} \frac{N(s)}{s-z} \frac{ds}{2\pi i} + \int^2_{-2} \frac{F(s)
 \delta^{\sigma_3}_{t+}(s) \left(\matrix 1&0\\ \ovl{r_t(s)}&1\endmatrix\right)}{s-z} \frac{ds}{2\pi i}.\tag\ppp{3.52}$$

Thus

$$\left\|M_-\delta^{\sigma_3}_{t-} \left(\matrix t&\bar r_t\\ 0&1\endmatrix\right)\right\|_{L^p(-1,1)}
\le \left\|\ \oint\limits_{|s|=2} \frac{N(s)}{s-\lozenge}
\frac{ds}{2\pi i}\right\|_{L^p(-1,1)}  +
c\|\delta^{\sigma_3}_{t+}\|_{L^\infty(\bb R)} \|F\|_{L^p(\bb
R)}.\tag\ppp{3.53}$$

Now from (\pp{3.46}), (\pp{3.47}), (\pp{3.49}), we have $M_- =
\text{I}' + \text{II}'$, where

$$\|\text{I}'\|_{L^p(\bb R)} \le \frac{c}{(1-\rho)^\beta} \|F\|_{L^p(\bb R)},
\|\text{II}''\|_{L^p(\bb R)} \le \frac{c}{(1-\rho)^{1+5\beta}}
\|F\|_{L^p(\bb R)},\tag\ppp{3.54}$$ and using (\pp{3.34}) we also
obtain $M_+ = \text{I}^0 + \text{II}^0$, where

$$\|\text{I}^0\|_{L^p(\bb R)} \le \frac{c}{(1-\rho)^\beta} \|F\|_{L^p(\bb R)}, \|\text{II}^0\|_{L^2(\bb R)}
 \le \frac{c}{(1-\rho)^{1+5\beta}} \|F\|_{L^p(\bb R)}.\tag\ppp{3.55}$$
But then using the Cauchy formula $M(z) = (C(M_+ - M_-))(z)$, we
can write $N$ as a sum of two parts, $\text{I}^N + \text{II}^N$,
and

$$\align
\|\text{I}^N\|_{L^p(|z|=2)} &\le \frac{c}{(1-\rho)^{\beta+\frac12}} \|F\|_{L^p(\bb R)},\\
\|\text{II}^N\|_{L^2(|z|=2)} &\le \frac{c}{(1-\rho)^{5\beta+3/2}}
\|F\|_{L^p(\bb R)}.
\endalign$$
The extra factor 1/2 comes from (\pp{3.13}), whereas $\bar r_t$
and $r_t$ are uniformly bounded in $\bb C_+$ and $\bb C_-$
respectively, for all $t>0$. Inserting these bounds in (\pp{3.53})
we obtain

$$\align
\|M_-\|_{L^p(-1,1)} &\le \frac{c}{(1-\rho)^{1/2}} \left(\frac1{(1-\rho)^{1/2}} + \frac1{(1-\rho)^{3/2+5\beta}}\right) \|F\|_{L^p}\\
&\le \frac{c}{(1-\rho)^{2+5\beta}} \|F\|_{L^p}.
\endalign$$
Together with (\pp{3.50}), this implies

$$\|M_-\|_{L^p(\bb R)} \le \frac{c}{(1-\rho)^{2+5\beta}} \|F\|_{L^p}\tag\ppp{3.56}$$
as $1+6\beta < 2+5\beta$ for $\beta < \frac12 < 1$. As before, the
same bound is true for all $F\in L^p(\bb R)$, by density. Finally,
by (\pp{1.13}) and (\pp{3.33}), this completes the proof of
Proposition (\pp{3.23}).

We now consider the general case where $r\in H^1_1$, $\|r\|_{H^1} \le \lambda$, $\|r\|_\infty \le \rho < 1$. We also continue to assume that $x=0$, and hence $z_0 = 0$ and $\theta = -tz^2$.

Given $r$ in $H^1_1$ as above, define

$$r^\#(z) = \frac{r(0)}{1+iz}.\tag\ppp{3.57}$$
Then $r^\#$ corresponds to a model problem of the form (\pp{3.22})
(ii) with $|r^\#(0)| = |r(0)| \le \rho < 1$, for which
Proposition~\pp{3.23} applies. Let $\delta,\delta^\#$ be the
solution of the scalar normalized RHP's (\pp{3.8}) associated with
$r,r^\#$ respectively. Set

$$\delta_1 = \delta(\delta^\#)^{-1}.\tag\ppp{3.58}$$
Then $\delta_1$ solves the scalar normalized RHP with jump

$$\delta_{1+} = \delta_{1-} \frac{1-|r|^2}{1-|r^\#|^2} \text{ on } \bb R_-.\tag\ppp{3.59}$$

By (\pp{3.19}), (\pp{3.20}), (\pp{3.21}), the jump matrix $\breve
v^\#_\theta$ associated with $r^\#$ has the form

$$\align
\breve v^\#_\theta &= \left(\matrix 1&r^\# e^{i\theta} \delta^{\#2}\\ 0&1\endmatrix\right)
\left(\matrix 1&0\\ -\bar r^\# e^{-i\theta}\delta^{\#-2}&1\endmatrix\right),\qquad z>0,\tag\ppp{3.60}\\
\breve v^\#_\theta &= \left(\matrix 1&0\\ \frac{-\bar r^\# e^{-i\theta}
\delta^{\#-2}}{1-|r^\#|^2}&1\endmatrix\right) \left(\matrix 1&\frac{r^\#e^{i\theta}
\delta^{\#2}}{1- |r^\#|^2}\\ 0&1\endmatrix\right)\tag\ppp{3.61}\\
&= \left(\matrix 1&0\\ -\bar r^\# e^{-i\theta} \delta^{\#-1}_+
\delta^{\#-1}_-&1\endmatrix\right) \left(\matrix 1&r^\#e^{i\theta}
\delta^\#_+ \delta^\#_-\\ 0&1\endmatrix\right),\qquad z<0.
\endalign$$
Set $v_1 = \delta^{\sigma_3}_{1-} \breve v^\#_\theta
\delta^{-\sigma_3}_{1+}$ and consider the IRHP$_{L^p}$

$$M_+  = M_-v_1 + F(v_1-I),\qquad F\in L^p(\bb R).\tag\ppp{3.62}$$
Using (\pp{1.14}), we obtain

$$M_- = [(1-C_{\breve v^\#_\theta})^{-1} C^-(F(\delta^{\sigma_3}_{1-} \breve v^\#_\theta -
\delta^{\sigma_3}_{1+}))] \delta^{-\sigma_3}_{1-}.$$
Then by (\pp{3.24}) and Proposition \pp{3.9},

$$\|M_-\|_{L^p} \le c\|(1-C_{\breve v^\#_\theta})^{-1} \|_{L^p\to L^p} \|F\|_{L^p} (\|\delta^{\sigma_3}_{1-}\|_{L^\infty} + \|\delta^{\sigma_3}_{1+}\|_{L^\infty}) \|\delta^{-\sigma_3}_{1-}\|_{L^\infty} \le \frac{c}{(1-\rho)^{3+5\beta}} \|F\|_{L^p},$$
and hence by (\pp{1.13})

$$\|(1-C_{v_1})^{-1}\|_{L^p\to L^p} \le \frac{c}{(1-\rho)^{3+5\beta}},\qquad 2<p<\infty.\tag\ppp{3.63}$$
Furthermore $v_1$ has the form $\left(\matrix 1-|r^\#|^2&
r^\#e^{i\theta}\delta^2\\ -\bar r^\#
e^{-i\theta}\delta^{-2}&1\endmatrix\right)$ for $z>0$ and
$\left(\matrix \frac{1-|r^\#|^2}{1-|r|^2}&r^\#
\delta_+\delta_-e^{i\theta}\\ -\bar r^\# \delta^{-1}_-
\delta^{-1}_+ e^{-i\theta}&1-|r|^2\endmatrix\right)$ for $z<0$.
Hence by (the proof of) Proposition \pp{3.2},

$$\|(1-C_{v_1})^{-1}\|_{L^2\to L^2} \le \frac{c}{1-\rho}.\tag\ppp{3.64}$$
Note that for $z>0$, $v_1$ can be written in the form

$$v_1 = \left(\matrix g&g(r^\# - r)e^{i\theta}\delta^2\\ 0&g^{-1}\endmatrix\right)
 \breve v_\theta \left(\matrix g&0\\ -g(\ovl{r^\#}-\bar r)e^{-i\theta} \delta^{-2}&g^{-1}\endmatrix\right)\tag\ppp{3.65}$$
and for $z<0$

$$v_1 = \left(\matrix g&0\\ -g^{-1}(\ovl{r^\#}-\bar r) e^{-i\theta} \delta^{-1}_+
\delta^{-1}_-&g^{-1}\endmatrix\right) \breve v_\theta
\left(\matrix g&g^{-1}(r^\#-r) e^{i\theta}\delta_+\delta_-\\
0&g^{-1} \endmatrix\right)\tag\ppp{3.66}$$ where

$$\align
g(z) &= \left(\frac{1-|r^\#(z)|^2}{1-|r(z)|^2}\right)^{\frac12} \quad \text{for}\quad z<0,\tag\ppp{3.67}\\
&= 1\quad \text{for}\quad z>0,
\endalign$$
and $\breve v_\theta = \delta^{\sigma_3}_- v_\theta
\delta^{-\sigma_3}_+$ as in (\pp{3.7}).

Extend $v_1$ to the complete, oriented contour $\Gamma = \bb R
\cup i\bb R$ on the RHS of Figure~\pp{1.35}, by setting

$$\align
v^e(z) &\equiv v_1(z),\qquad z>0\tag\ppp{3.68}\\
&\equiv v^{-1}_1(z),\qquad z<0\\
&\equiv I,\qquad z\in i\bb R.
\endalign$$

The following fact is simple to prove (cf.\ [DZ4], [DZ5]). Suppose
$\Sigma$ is an oriented contour in $\bb C$ with associated jump
matrix $v$, and suppose we reverse the orientation on some subset
$\Sigma'\subset \Sigma$. Denote the new contour by
$\widehat\Sigma$ and set $\hat v = v$ on $\Sigma\backslash
\Sigma'$, $\hat v = v^{-1}$ on $\Sigma'$. Then the operators
$C_{\hat v}$ and $C_v$ on $L^p(\widehat\Sigma)$ and $L^p(\Sigma)$
respectively, $1<p<\infty$, are the same i.e.\ $C_{\hat v}f = C_v
f$ for all $f\in L^p(\widehat \Sigma) \equiv L^p(\Sigma)$.
Together with the fact that $v^e\equiv I$ on $\Gamma \backslash
\bb R$, this implies by (\pp{3.63}), (\pp{3.64}) that

$$\|(1-C_{v^e})^{-1}\|_{L^p(\Gamma)} \le \frac{c}{(1-\rho)^{3+5\beta}},\qquad 2<p<\infty,\tag\ppp{3.69}$$
and

$$\|(1-C_{v^e})^{-1}\|_{L^2(\Gamma)} \le \frac{c}{1-\rho}\tag\ppp{3.70}$$
for some constants $c$.

Let $\Omega_j = \{z\colon\ (j-1) \frac\pi2 < \arg z < j\frac\pi2\}$, $1\le j\le 4$, denote the four components of $\bb C\backslash \Gamma$ with oriented boundaries

$$\align
\Sigma_1 &= \{+i\infty \to 0 \to +\infty\}\\
\Sigma_2 &= \{+i\infty \to 0 \to -\infty\}\\
\Sigma_3 &= \{-i\infty\to 0 \to -\infty\}\\
\Sigma_4 &= \{-i\infty \to  0 \to +\infty\}
\endalign$$
respectively. Note that in the notation of Figure \pp{1.35},
$\Omega_+ = \Omega_1 \cup \Omega_3$ and $\Omega_-  = \Omega_2 \cup
\Omega_4$. Set

$$\align
v_2(z) &= \breve v_\theta(z),\qquad z>0,\tag\ppp{3.71}\\
&= \breve v^{-1}_\theta(z),\qquad z<0,\\
&= I,\qquad z\in i\bb R.
\endalign$$
With this notation, we have

$$v^e = G_- v_2 G_+\tag\ppp{3.72}$$
where

$$\align
G_+ &= G_1 \quad \text{on}\quad \Sigma_1,\tag\ppp{3.73}\\
G_1(z) &= \left(\matrix 1&0\\ 0&1\endmatrix\right),\qquad z\in i\bb R_+,\\
&= \left(\matrix g&0\\ -g(\ovl{r^\#}-\bar r)e^{-i\theta}
\delta^{-2}&g^{-1}\endmatrix\right),\qquad z>0,
\endalign$$

$$\align
G_+ &= G_3 \quad \text{on}\quad \Sigma_3,\tag\ppp{3.74}\\
G_3(z) &= \left(\matrix 1&0\\ 0&1\endmatrix\right), \qquad z\in i\bb R_-,\\
&= \left(\matrix g^{-1}&0\\ -g^{-1}(\ovl{r^\#}-\bar r)e^{-i\theta}
\delta^{-1}_+\delta^{-1}_-&g\endmatrix\right), \qquad z<0,
\endalign$$

$$\align
G_- &= G_2 \quad \text{on}\quad \Sigma_2,\tag\ppp{3.75}\\
G_2(z) &= \left(\matrix 1&0\\ 0&1\endmatrix\right), \qquad z\in i\bb R_+,\\
&= \left(\matrix g^{-1}&-g^{-1}(r^\#-r) e^{i\theta}
\delta_+\delta_-\\ 0&g\endmatrix\right),\qquad z<0,
\endalign$$
and

$$\align
G_- &= G_4 \quad \text{on}\quad \Sigma_4,\tag\ppp{3.76}\\
G_4(z) &= \left(\matrix 1&0\\ 0&1\endmatrix\right) \quad \text{on}\quad i\bb R_-\\
&= \left(\matrix g&g(r^\# -r)e^{i\theta} \delta^2\\
0&g^{-1}\endmatrix\right),\qquad z>0.
\endalign$$
Note that $\|G_\pm\|_{L^\infty(\Gamma)} \le
\frac{c}{(1-\rho)^{1/2}}$.

Now extend $\Gamma$ to the complete, oriented contour
$\Gamma_{\text{ext}} = \Gamma\cup \Gamma'$ on the RHS of
Figure~\pp{1.38}. Choosing $\text{dist}(\Gamma,\Gamma')<1$
sufficiently small, we will show how to construct an invertible
matrix function $H$ with $H$ and $H^{-1}$ analytic and bounded in
$\bb C\backslash \Gamma_{\text{ext}}$, such that $H_\pm-G_\pm$ is
as small as desired on $\Gamma$. Here $H_\pm$ denote the boundary
values of $H(z)$ on $\Gamma$ from $\Omega_\pm$ respectively (cf.\
Figure~\pp{1.38}). It then follows by perturbation theory that the
bounds (\pp{3.69}), (\pp{3.70}) remain true for the Cauchy
operator with $v^e = G_-v_2G_+$ replaced by $v^H \equiv
H_-v_2H_+$. But then as $H$ is piecewise analytic, it follows by
the Conjugation Lemma~\pp{1.39} that similar bounds are true for
$v_2$, and hence for $\breve v_\theta$, and hence for $v_\theta$,
as desired.

We will show how to construct $H = H(z)$ in $\Omega_2$. The construction of $H$ in $\Omega_1, \Omega_3$ and $\Omega_4$ is similar and left to the reader. As above, $r\in H^1_1$, $\|r\|_{H^1} \le \lambda$, $\|r\|_{L^\infty} \le \rho < 1$.

Set $g(z) =1$ on $i\bb R_+$. Note that with this definition (cf.\
(\pp{3.67})) $g(z)$ is continuous on $\Sigma_2$. Note also that
the same is true for the function which equals $-g^{-1}(r^\#-r)
\delta_+\delta_-$ on $\bb R_-$ and equals zero on $i\bb R_+$. For
some constant $c$, any $0<\vp <1$, and $\gamma>0$ satisfying

$$\sqrt\gamma < \frac{\vp(1-\rho)^{5/2}}{c\lambda},\tag\ppp{3.77}$$
we construct functions $H_2,h_2$ with the following properties:

\medskip

\n (\ppp{3.78})

\itemitem{(i)}
$\|H^{-\sigma_3}_2-g^{-\sigma_3}\|_{L^\infty(\Sigma_2)} < \vp$

\itemitem{(ii)} $H_2(z)\to 1\quad \text{as}\quad z\to \infty \quad \text{in}\quad \Sigma_2$

\itemitem{(iii)} $\|H^{\sigma_3}_2-I\|_{L^2(\Sigma_2)} \le \frac{c\lambda}{1-\rho}$

\itemitem{(iv)} $H^{\sigma_3}_2$ is analytic in $\Sigma_{2,\gamma}
=\{z\colon \ \text{dist}(z,\Sigma_2) < \frac\gamma2\}$ and
$\|H^{\sigma_3}_2\|_{L^\infty(\Sigma_{2,\gamma})} \le
\frac{c}{(1-\rho)^{1/2}}$

\itemitem{(i)$'$} $\|h_2\|_{L^\infty(i\bb R_+)} + \|h_2 - (-g^{-1}(r^\# -r)\delta_+\delta_-)\|_{L^\infty(\bb R_-)} < \vp$

\itemitem{(ii)$'$} $h_2(z) \to 0$ as $z\to\infty$ in $\Sigma_2$

\itemitem{(iii)$'$} $\|h_2\|_{L^2(\Sigma_2)} \le \frac{c\lambda}{(1-\rho)^{1/2}}$

\itemitem{(iv)$'$} $h_2$ is analytic in $\Sigma_{2,\gamma}$ and $\|h_2\|_{L^\infty(\Sigma_{2,\gamma})} \le
 \frac{c}{(1-\rho)^{1/2}}$.

Define $H_2$ as follows. For $z\in \Sigma_2$,

$$H_2(z) \equiv - \int_{\Sigma_2} g(s) \left(\frac1{s-(z + \gamma e^{3i\pi/4})} -
\frac1{s-(z-\gamma e^{3i\pi/4})}\right) \frac{ds}{2\pi
i}.\tag\ppp{3.79}$$ Using the fact that $1 = -\int_{\Sigma_2}
\left(\frac1{s-(z+\gamma e^{i3\pi/4})} - \frac1{s-(z-\gamma
e^{i3\pi/4})}\right) \frac{ds}{2\pi i}$, we obtain

$$|H_2(z) - g(z)| \le \int_{\Sigma_2} \frac{|g(s) - g(z)|}{|(s-z)^2 + i\gamma^2|} \frac{\gamma|ds|}{\pi},\tag\ppp{3.80}$$
and inserting the estimate $|g(s) - g(z)| \le |s-z|^{1/2}
\|g'\|_{L^2(\bb R_-)}$ for $s,z\in \Sigma_2$, we arrive at the
estimate $|H_2(z) - g(z)| \le c\sqrt\gamma\|g'\|_{L^2(\bb R_-)}$.
A simple calculation yields $\|g'\|_{L^2(\bb R_-)} \le
\frac\lambda{(1-\rho)^{3/2}}$ and hence

$$\|H_2-g\|_{L^\infty(\Sigma_2)} \le \frac{c\sqrt \gamma\ \lambda}{(1-\rho)^{3/2}}.\tag\ppp{3.81}$$
Using the fact that $|g(z)| \ge \sqrt{1-\rho}$ on $\Sigma_2$, a
straightforward calculation shows that if we choose $\gamma>0$
such that

$$\frac{2\sqrt 2 c \sqrt \gamma \ \lambda}{(1-\rho)^{5/2}} < \vp < 1,\tag\ppp{3.82}$$
then $|H_2(z)| > \frac12\sqrt{1-\rho}$ and

$$\|H^{\sigma_3}_2 - g^{\sigma_3}\|_{L^\infty(\Sigma_2)} \le \frac\vp{\sqrt 2} \|\text{I}\|  = \vp\tag\ppp{3.83}$$
(recall Remark \pp{1.56}), which proves (i). The proof of (ii) is
elementary and replacing $g(z)$ by 1 in (\pp{3.80}) we obtain by a
standard calculation the bound $\|H_2-1\|_{L^2(\Sigma_2)} \le
c\|g-1\|_{L^2(\Sigma_2)}$. But then direct estimation shows that
$\|g-1\|_{L^2(\Sigma_2)} \le \lambda/\sqrt{1-\rho}$ and so using
the above fact that $|H_2(z)| \ge \frac12\sqrt{1-\rho}$, we obtain
(iii),

$$\|H^{\sigma_3}_2 -I\|_{L^2(\Sigma_2)} \le \frac{c\lambda}{1-\rho}.\tag\ppp{3.84}$$
Finally, it is clear from (\pp{3.79}) that $H_2(z)$ extends to an
analytic function in $\Sigma_{2,\gamma}$, and inserting $z+\mu
e^{i3\pi/2}$, $-\frac\gamma2 < \mu  < \frac\gamma2$, we obtain as
before the bound $|H_2(z+\mu e^{i3\pi/4}) - g(z)| \le
c\sqrt\gamma\|g'\|_{L^2(\bb R_-)} \le \frac{c\sqrt \gamma\
\lambda}{(1-\rho)^{3/2}}$ for $z\in\Sigma_2$. Using (\pp{3.83})
this then leads to (iv),

$$\|H^{\sigma_3}_2\|_{L^\infty(\Sigma_{2,\gamma})} \le \frac{c}{(1-\rho)^{1/2}}\tag\ppp{3.85}$$
as desired.

We now construct $h_2$. Define $b$ on $\Sigma_2$ by $b(s) = 0$ for
$s\in i\bb R_+$ and $b(s) = (r-r^\#)g^{-1} \Delta(s)$ for $s\in\bb
R_-$, where $\Delta = \delta_+\delta_-$ as in (\pp{3.17}). For
$z\in \Sigma_2$, set

$$h_2(z) = -\int_{\Sigma_2} b(s) \left(\frac1{s-(z + \gamma e^{i3\pi/4})} - \frac1{s-(z-\gamma e^{i3\pi/4})}\right)
\frac{ds}{2\pi i}.\tag\ppp{3.86}$$ As above, we obtain
$\|h_2-b\|_{L^\infty(\Sigma_2)} \le c\sqrt\gamma\|b'\|_{L^2(\bb
R_-)}$. Now for $z<0$, $b' = (r'-(r^\#)')g^{-1} \Delta +\break
(r-r^\#) (g^{-1})'\Delta + (r-r^\#)g^{-1} \Delta' = \text{I} +
\text{II} + \text{III}$. Clearly $\|\text{I}\|_{L^2(\bb R_-)} \le
\frac{c\lambda}{\sqrt{1-\rho}}$. As $g$ and $g^{-1}$ have the same
structure we obtain as above $\|(g^{-1})'\|_{L^2(\bb R_-)} \le
\frac\lambda{(1-\rho)^{3/2}}$ and so $\|\text{II}\|_{L^2(\bb R_-)}
\le\frac{c\lambda}{(1-\rho)^{3/2}}$. Now from (\pp{3.10}),
(\pp{3.17}), $\Delta = e^{-H((\log(1-|r|^2))\chi_{\bb R_-})}$,
where $H = -(C^++C^-)$ is the Hilbert transform and $\chi_{\bb
R_-}$ is the characteristic function for $\bb R_-$. Hence

$$\Delta'(z) = -\Delta(z) \frac{d}{dz} H((\log(1-|r|^2)) \chi_{\bb R_-}) = \Delta(z) H\left(\frac{|r|^{2'}}{1-|r|^2}
 \chi_{\bb R_-}\right) -\frac{i\Delta}{\pi} \frac{\log (1-|r(0)|^2)}{z}.$$
Using the $L^2$ mapping properties of $H$ (see (\pp{1.1}) (ii)),
the identity $|\Delta| = 1$, and the elementary bound $|\log(1 -
|r(0)|^2)| \le \frac{|r(0)|^2}{1-|r(0)|^2}$, we see that
$|\Delta'(z)|\le \text{I}' + \text{II}'$ where
$\|\text{I}'\|_{L^2} \le \frac{c\lambda}{1-\rho}$ and
$|\text{II}'(z)| \le \frac{c}{1-\rho} \frac1{|z|}$, $z<0$. Thus
$\|\text{III}\|_{L^2(\bb R_-)} \le \frac{c\lambda}{(1-\rho)^{3/2}}
+ \frac{c}{\sqrt{1-\rho}} \cdot \frac1{1-\rho}
\left\|\frac{r-r^\#}{|\cdot |}\right\|_{L^2(\bb R)}$. But then by
Hardy's inequality, $\left\|\frac{r-r^\#}{|\cdot
|}\right\|_{L^2(\bb R_-)} \le\break 2 \|r'-(r^\#)'\|_{L^2(\bb
R_-)} \le 4\lambda$, and so $\|\text{III}\|_{L^2(\bb R_-)} \le
\frac{c\lambda}{(1-\rho)^{3/2}}$, and hence $\|b'\|_{L^2(\bb R_-)}
\le \frac{c\lambda}{(1-\rho)^{3/2}}$. Choosing
$\frac{c\sqrt\gamma\ \lambda}{(1-\rho)^{3/2}} <\vp$, we obtain
(i)$'$. The proof of (ii)$'$ is immediate and as
$\|h_2\|_{L^2(\Sigma_2)} \le c\|b\|_{L^2(\Sigma_2)} \le
\frac{c\lambda}{\sqrt{1-\rho}}$, we obtain (iii)$'$. Again it is
clear from (3.86) that $h_2(z)$ extends to an analytic function in
$\Sigma_{2,\gamma}$, and inserting $z + \mu e^{i3\pi/4}$,
$-\frac\gamma2 < \mu < \frac\gamma2$, we see that $|h_2(z + \mu
e^{i3\pi/4})| \le c \|b\|_{L^\infty(\Sigma_2)} \le
\frac{c}{\sqrt{1-\rho}}$. This proves (iv)$'$. Adjusting the
constants $c$ at various points in the above construction, we see
that for $0<\vp<1$ and $\gamma$ as in (\pp{3.77}) we have obtained
$H_2$, $h_2$ with the desired properties (\pp{3.78}).

Define $H$ as a piecewise analytic matrix function in $\Omega_2$
as follows:\ For $z\in \Omega_2 \cap \Sigma_ {2,\gamma} \subset
\Omega_-$ (see cf.\ Figure~\pp{1.38}) we set

$$H(z) = \left(\matrix H^{-1}_2(z)&h_2(z)e^{i\theta}\\ 0&H_2(z)\endmatrix\right)\tag\ppp{3.87}$$
and for $z\in \Omega_2\backslash \Sigma_{2\gamma} \subset
\Omega_+$ (cf.\ Figure \pp{1.38}) we set

$$H(z) = I.\tag\ppp{3.88}$$
Similar constructions taking into account the triangularity of
$G_\pm$ (see (\pp{3.73}), (\pp{3.74}), (\pp{3.76})) yield $H$ in
$\Omega_1$, $\Omega_3$ and $\Omega_4$ respectively:\ the details
are left to the reader. We obtain an invertible matrix valued
function $H$ on $\bb C\backslash \Gamma_{\text{ext}}$ with the
following properties for some constant $c>0$, and

$$\gather
0 < \frac{c\lambda\sqrt\gamma}{(1-\rho)^{5/2}} < \vp < 1,\tag\ppp{3.89}\\
H(z)\quad \text{is analytic and bounded in}\quad \bb C\backslash
\Gamma_{\text{ext}},\
\|H^{\pm1}\|_{\bb C\backslash \Gamma_{\text{ext}}} \le \frac{c}{(1-\rho)^{1/2}},\tag\ppp{3.90}\\
\|H_\pm - G_\pm\|_{L^\infty(\Gamma)} < \vp,\tag\ppp{3.91}
\endgather$$
where $H_\pm$ denote the boundary values of $H(z)$ on the oriented
contour $\Gamma\subset \Gamma_{\text{ext}}.$

$$\gather
H(z) \to I \text{ uniformly as } z\to \infty \text{ in } \bb C\backslash \Gamma_{\text{ext}} \text{ and}\tag\ppp{3.92}\\
\|H_\pm - I\|_{L^2(\Gamma)},\quad \|H^{-1}_\pm-I\|_{L^2(\Gamma)}
\le \frac{c\lambda}{1-\rho}.
\endgather$$
Note that in deriving these properties the signature table for
$\text{Re } i\theta$ in Figure~\pp{3.5} plays a crucial role.

Now if $\frac{c\sqrt \gamma\ \lambda}{(1-\rho)^{5/2}} < \vp$, then

$$\align
\|v^H-v^e\|_{L^\infty(\Gamma)} &= \|H_-v_2H_+ - G_-v_2G_+\|_{L^\infty(\Gamma)}\tag\ppp{3.93}\\
&\le \frac{c\vp}{\sqrt{1-\rho}} \|v_2\|_{L^\infty(\Gamma )} \le
\frac{c\vp}{\sqrt{1-\rho}}
\endalign$$
by (\pp{3.71}), (\pp{3.19}), and (\pp{3.21}). Thus for
$2<p<\infty$, by (\pp{3.69}),

$$\|(v^H -v^e)(1-C_{v^e})^{-1}\|_{L^p(\Gamma)\to L^p(\Gamma)} \le \frac{c\vp}{\sqrt{1-\rho}} \cdot \frac{c}{(1-\rho)^{3+5\beta}} < \frac12$$
if $\vp < c(1-\rho)^{7/2+5\beta}$. Thus if
$$\sqrt\gamma < \frac{c(1-\rho)^{5/2}}{\lambda} \cdot c(1-\rho)^{7/2+5\beta} =
\frac{c(1-\rho)^{6+5\beta}}{\lambda},
\tag\ppp{3.94}$$ we see by the second resolvent identity that
$(1-C_{v^H})^{-1}$ exists in $L^p(\Gamma)$ and

$$\|(1-C_{v^H})^{-1}\|_{L^p(\Gamma)\to L^p(\Gamma)} \le \frac{c}{(1-\rho)^{3+5\beta}}.\tag\ppp{3.95}$$
Similarly for $\gamma$ satisfying (\pp{3.94}) (adjust $c$ if
necessary) we have from (\pp{3.70})

$$\|(1-C_{v^H})^{-1}\|_{L^2(\Gamma)\to L^2(\Gamma)} \le \frac{c}{(1-\rho)}.\tag\ppp{3.96}$$

We are now in a position to apply the Conjugation Lemma \pp{1.39}.
In $\bb C\backslash \Gamma_{\text{ext}}$ (see Figure~\pp{1.38}),
set $R(z) = H^{-1}(z)$ for $z\in \Omega_+$ and $R(z) = H(z)$ for
$z\in \Omega_-$. Then $v_2 = R^{-1}_- v^H R_+$. Clearly $R(z)$ is
analytic and invertible in $\bb C\backslash \Gamma_{\text{ext}}$
and $\|R\|_{L^\infty(\bb C\backslash \Gamma_{\text{ext}})}$,
$\|R^{-1}\|_{L^\infty(\bb C\backslash \Gamma_{\text{ext}})} \le
\frac{c}{(1-\rho)^{1/2}}$. As in (\pp{3.43}), we have \break
$\|(1-C_{\breve v_\theta})^{-1}\|_{L^2(\bb R)} \le
\frac{c}{1-\rho}$, and taking into account the discussion
preceding (\pp{3.69}), (\pp{3.70}), we obtain
$\|(1-C_{v_2})^{-1}\|_{L^2(\Gamma)} \le \frac{c}{1-\rho}$. Also
$\|v_2\|_{L^\infty(\Gamma)}\le c$, $\|v^H\|_{L^\infty} \le
\frac{c}{1-\rho}$, $\|v_2-I\|_{L^2(\Gamma)} \le \lambda$,
$\|v^H-I\|_{L^2(\Gamma)} \le \frac{c\lambda}{(1-\rho)^{3/2}}$, and
provided (\pp{3.94}) holds, (\pp{3.95}), (\pp{3.96}) are
satisfied. It then follows by (\pp{1.44}) that

$$\|(1-C_{v_2})^{-1}\|_{L^p(\Gamma)\to L^p(\Gamma)} \le c^\#,\tag\ppp{3.97}$$
where

$$\align
c^\# &= c_{\text{dist}(\Gamma,\Gamma')} \|R\|_{L^\infty( \bb C\backslash\Gamma_{\text{ext}})}
\|R^{-1}\|_{L^\infty(\bb C\backslash \Gamma_{\text{ext}})}\|(1-C_{v^H})^{-1}\|_{L^p(\Gamma)} \|(1-C_{v^H})^{-1}\|_{L^2(\Gamma)}\tag\ppp{3.98}\\
&\quad \times \|(1-C_{v_2})^{-1}\|^2_{L^2(\Gamma)}
\|v_2\|^3_{L^\infty(\Gamma)} \|v^H\|_{L^\infty(\Gamma)}^2
 (1+ \|v^H-I\|_{L^2(\Gamma)})^2
(1 + \|v_2-I\|_{L^2(\Gamma)})^2\\
&\le c_{\text{dist}(\Gamma,\Gamma')} \frac1{(1-\rho)^{\frac12}} \frac1{(1-\rho)^{\frac12}}
\frac1{(1-\rho)^{3+5\beta}} \frac1{(1-\rho)} \frac1{(1-\rho)^2}
\frac1{(1-\rho)^2}\left(1+ \frac\lambda{(1-\rho)^{3/2}}\right)^2\\
&\quad \times (1+\lambda)^2\\
&\le c_{\text{dist}(\Gamma,\Gamma')}
\frac{(1+\lambda)^4}{(1-\rho)^{12+5\beta}},
\endalign$$
where $c_{\text{dist}(\Gamma,\Gamma')}$ depends on the distance
$\frac\gamma2$ between $\Gamma$ and $\Gamma'$. By (\pp{2.14}),
$c_{\text{dist}(\Gamma,\Gamma')} = \frac{c_p}{\gamma^{3/2+1/p}}$,
provided $\frac\gamma2 \le 1$. Choose

$$\sqrt\gamma = \frac{c(1-\rho)^{6+5\beta}}{1+\lambda}, \tag\ppp{3.99}$$
where the constant $c$ may be taken as the minimum of 1 and the
constant on the RHS of (\pp{3.94}). Then certainly $\frac\gamma2 <
\gamma <1$, and we obtain the bound

$$\align
\|(1-C_{v_2})^{-1}\|_{L^p\to L^p} &\le \frac{c_p
(1+\lambda)^4}{(1-\rho)^{12+5\beta}}
 \left(\frac{(1+\lambda)^2}{(1-\rho)^{12+10\beta}}\right)^{\frac32+\frac1p}\tag\ppp{3.100}\\
&=  \frac{c_p
(1+\lambda)^{7+\frac2p}}{(1-\rho)^{30+\frac{12}p+\beta'}},
\endalign$$
where $\beta'>0$. But then reversing the discussion preceding
(\pp{3.69}), (\pp{3.70}), we obtain for $2<p<\infty$,

$$\|(1-C_{\breve v_\theta})^{-1}\|_{L^p(\bb R)\to L^p(\rb)} \le \frac{c_p (1+\lambda)^{7+\frac2p}}
{(1-\rho)^{30+\frac{12}p + \beta'}},\qquad
\beta'>0.\tag\ppp{3.101}$$ Using the formula

$$(1-C_{v_\theta})^{-1}h = [(1-C_{\breve v_\theta})^{-1} (C^-(h(v_\theta-I) \delta^{-\sigma_3}_+))]
\delta^{\sigma_3}_- + h$$
which follows from (\pp{1.13}), (\pp{1.14}), we obtain finally
that for $2<p<\infty$

$$\|(1-C_{v_\theta})^{-1}\|_{L^p(\bb R)\to L^p(\rb)} \le \frac{c_p(1+\lambda)^{7+\frac2p}}
{(1-\rho)^{31+\frac{12}p + \beta'}},\qquad
\beta'>0.\tag\ppp{3.102}$$ By (\pp{3.3}), we have

$$\|(1-C_{v_\theta})^{-1}\|_{L^2(\bb R)\to L^2(\rb)} \le \frac{c}{1-\rho}.\tag\ppp{3.103}$$

This completes the proof of (\pp{1.49}) in the case that $x=0$, for
all $t\in\bb R$. For general $x\ne 0$, and $z_0 = \frac{x}{2t}$,
let $T_{z_0}f(z) = f(z+z_0)$. Then $T_{z_0}$ is an isometry from
$L^p(\bb R)$ onto $L^p(\bb R)$ and a direct calculation shows that
$T_{z_0}\circ C_{v_\theta} \circ T^{-1}_{z_0} =
C_{v_{z_0,\theta}}$ where

$$v_{z_0,\theta}(z)  = \left(\matrix 1-|r_{z_0}(z)|^2&r_{z_0}(z)e^{-itz^2}\\
-\ovl{r_{z_0}(z)} e^{itz^2}&1\endmatrix\right)$$ and $r_{z_0}(z) =
r(z+z_0) e^{itz_0^2}$. As $r_{z_0} \in H^1_1(\bb R)$ and
$\|r_{z_0}\|_{H^1(\bb R)} = \|r\|_{H^1} \le \lambda$,
$\|r_{z_0}\|_{L^\infty(\bb R)} = \|r\|_{L^\infty(\bb R)} \le\rho
<1$, the general case now follows from the case $x=0$. Thus
(\pp{3.102}), (\pp{3.103}) are true for all $x,t\in \bb R$.

We can apply Riesz--Thorin interpolation to (\pp{3.102}),
(\pp{3.103}). For any $k>1$, $2<p<\infty$, we find (denote the
constant in the $L^2$ bound by $c_2$)

$$\|(1-C_{v_\theta})^{-1}\|_{L^p(\bb R)\to L^p(\rb)} \le \frac{c^{\xi}_{kp} c^{1-\xi}_2 (1+\lambda)^{(7+\frac2p)\xi}}{(
1-\rho)^{(31+\frac{12}p + \beta')\xi+(1-\xi)}}$$ where $\xi =
\left(1-\frac2p\right) \left(1-\frac2{kp}\right)$. For example,
for $p=4$ which is a case of principal interest in [DZ5], given
$\nu>0$, we can choose $\beta'$ sufficiently small and $k$
sufficiently large so that

$$\|(1-C_{v_\theta})^{-1}\|_{L^4(\bb R)\to L^4(\rb)} \le \frac{c'_4(1+\lambda)^{\frac{15}4 +\nu}}{(1-\rho)^{\frac{35}2+\nu}}\tag\ppp{3.104}$$
for some (very large) constant $c'_4$, which should be compared
with (\pp{3.102}) for $p=4$.

\remark{Remark}

In [DZ5], one needs a bound of the form
$\|(1-C_{w_\theta})^{-1}\|_{L^p(\bb R)\to L^p(\rb)} \le c$ for all
$x\in\bb R$ and for all $t\ge t_0$, where the time $t_0$ is large.
This is a much simpler situation than considered in this paper:\
to prove this bound one still uses steepest descent methods, but
the Conjugation Lemma~\pp{1.39} is not needed. We refer the reader
to [DZ5] for the details.

\endremark

\enddocument